\documentclass[a4paper, colorinlistoftodos]{article}

\usepackage[english]{babel}
\usepackage{csquotes}
\DeclareUnicodeCharacter{0393}{$\Gamma$}
\DeclareUnicodeCharacter{2081}{$_1$}
\DeclareUnicodeCharacter{2113}{$\ell$}

\usepackage{amsmath}
\usepackage{mathrsfs}
\usepackage{mathtools}
\usepackage{multicol}
\usepackage{amssymb}
\usepackage{enumerate}
\usepackage{dsfont}
\usepackage{bm}
\usepackage{graphicx}
\usepackage{hyperref}
\usepackage{physics}
\usepackage{siunitx}
\usepackage{multicol}
\usepackage{xcolor}
\usepackage{centernot}
\usepackage{appendix}
\usepackage{authblk}
\usepackage[a4paper, lmargin=2.5cm, rmargin=2.5cm, marginpar=3.0cm]{geometry}

\DeclarePairedDelimiterX\Set[1]\lbrace\rbrace{\def\given{\;\delimsize\vert\;\allowbreak}#1}

\usepackage[linesnumbered, ruled]{algorithm2e}
\DontPrintSemicolon
\SetKwComment{Comment}{/* }{ */}

\usepackage{amsthm}

\usepackage{hyperref}
\usepackage{cleveref}
\usepackage{breakurl}

\newcommand{\R}{\mathbb{R}}

\DeclareMathOperator*{\argmin}{arg\,min}
\DeclareMathOperator*{\sign}{sgn}
\DeclareMathOperator*{\diag}{diag}
\DeclareMathOperator*{\prox}{prox}
\DeclareMathOperator*{\lip}{Lip}

\renewcommand{\i}{\mathrm{i}}
\newcommand{\e}{\mathrm{e}}

\newcommand{\A}{\mathbf{A}}
\newcommand{\B}{\mathbf{B}}
\newcommand{\W}{\mathbf{W}}
\renewcommand{\b}{\mathbf{b}}

\newcommand{\X}{\mathcal{X}}
\newcommand{\Y}{\mathcal{Y}}
\newcommand{\U}{\mathcal{U}}

\newcommand{\RC}{\mathcal{R}}

\newcommand{\phiAE}{\bm{\phi}_{\theta}}
\newcommand{\phiEnc}{\bm{\phi}_{\theta,\text{enc}}}
\newcommand{\phiDec}{\bm{\phi}_{\theta,\text{dec}}}

\usepackage{caption}
\usepackage{subcaption}
\usepackage{svg}

\usepackage{arydshln}
\usepackage{lastpage}
\usepackage{fancyhdr}

\allowdisplaybreaks

\usepackage{todonotes}
\numberwithin{equation}{section}

\usepackage{parskip}

\usepackage[
    backend=biber,natbib,
    style=apa
]{biblatex}

\usepackage{xpatch}
  \xpretobibmacro{doi+url}{\printfield{isbn}}{}{}

\title{Sparsifying dimensionality reduction of PDE solution data\\ with Bregman learning}
\author[1,*]{Tjeerd Jan Heeringa}
\author[1]{Christoph Brune}
\author[1,2]{Mengwu Guo}
\affil[1]{Mathematics of Imaging \& AI, University of Twente, Enschede, The Netherlands}
\affil[2]{Centre for Mathematical Sciences, Lund University, Lund, Sweden}
\affil[*]{Corresponding author: t.j.heeringa@utwente.nl}
\date{}

\begin{document}
\maketitle

\begin{abstract}
\noindent Classical model reduction techniques project the governing equations onto a linear subspace of the original state space. More recent data-driven techniques use neural networks to enable nonlinear projections. Whilst those often enable stronger compression, they may have redundant parameters and lead to suboptimal latent dimensionality. To overcome these, we propose a multistep algorithm that induces sparsity in the encoder-decoder networks for effective reduction in the number of parameters and additional compression of the latent space. This algorithm starts with sparsely initialized a network and training it using linearized Bregman iterations. These iterations have been very successful in computer vision and compressed sensing tasks, but have not yet been used for reduced-order modelling. After the training, we further compress the latent space dimensionality by using a form of proper orthogonal decomposition. Last, we use a bias propagation technique to change the induced sparsity into an effective reduction of parameters. We apply this algorithm to three representative PDE models: 1D diffusion, 1D advection, and 2D reaction-diffusion. Compared to conventional training methods like Adam, the proposed method achieves similar accuracy with 30\% less parameters and a significantly smaller latent space. 
\\ \\
\textbf{keywords: }Nonlinear dimensionality reduction, linearized Bregman iterations, sparsity, neural architecture search, scientific machine learning

\end{abstract}

\section{Introduction}
Physically consistent simulation of processes is an essential tool used in science and engineering. To get high-fidelity simulations it is often required to use a fine spatio-temporal discretization; this can lead to very large models which take very long to solve. However, in many cases we want a real-time solution or solve the model for several different parameters. 

One of the ways to overcome this computational burden is to use projection-based reduced-order models. These techniques use a two-stage approach. In the \textit{offline} stage, a representative low-dimensional space is sought. This typically involves solving the full-order model for several parameters. These are expensive computations, but only have to be executed once. In the \textit{online} stage, solutions are computed in the latent space. Finding a reduced-order latent solution can be significantly faster than finding a full-order solution. The latent solutions are subsequently used to construct the full-order solutions. 

A lot of the projection-based use linear maps to describe an encoding-decoding structure. The best-known example of this is the \textit{Proper Orthogonal Decomposition} (POD) [\cite{sirovich_turbulence_1987,quarteroni_reduced_2016}], also known as Principal Component Analysis in statistical analysis, Empirical Orthogonal Eigenfunctions in atmospheric modelling and as Karkunen-Loève expansion in stochastic process modelling [\cite{hotelling_analysis_1933,loeve_probability_1963,north_sampling_1982}]. It uses the singular value decomposition of a data matrix to construct the linear maps and controls the error made by the reduction through the decay of the singular values. Several extensions have been developed to cover cases like the data matrix being too large for memory or containing missing/corrupted data [\cite{prudhomme_reliable_2001, everson_karhunenloeve_1995}]. Other methods for linear projection include Multidimensional Scaling (MDS), Linear Discriminant Analysis (LDS), Canonical Correlation Analysis (CCA), Maximum Autocorrection Factors (MAF), Sufficient Dimension Reduction (SDR), Locality Preserving Projections, and Least Squares [\cite{cunningham_linear_2015}]. For problems with a Kolgomorov $n$-width that decays fast enough in $n$ the usage of these approaches is suitable. However, many problems of interest have a Kolgomorov $n$-width that decays too slow [\cite{lee_model_2020,peherstorfer_breaking_2022}]. In these problems, using linear projection-based approaches provide limited value when the full-order solution has a high dimensionality.

To overcome this limitation set by the Kolgomorov $n$-width, various nonlinear methods have been developed. Well-known examples include Dynamic Mode Decomposition (DMD) [\cite{schmid_dynamic_2008,tu_dynamic_2014}] and Kernel Proper Orthogonal Decomposition (kPOD) [\cite{scholkopf_nonlinear_1998}]. In recent years, neural networks have become one of the main tools for creating new reduced-order methods [\cite{carlberg_galerkin_2017,fresca_comprehensive_2021}] due to their universal approximation property. This allows them to be used as replacement for any function in an algorithm, e.g. an autoencoder is a generalization of POD where the matrices used for reduction and recovery are replaced by neural networks. These function approximations tend to become more accurate when more data are available and when the number of their parameters is increased. Unfortunately, the memory requirements for storing the networks can become quite large and their general and nonlinear nature means that their behaviour is often hard to explain [\cite{nauta_anecdotal_2022}].

\subsection{Related work on sparse models}
Parsimonious models allow for faster and less memory-intensive reduction and recovery maps and smaller latent spaces. The lottery ticket hypothesis tells us that each network has a subnetwork that performs just as well as the original network [\cite{frankle_lottery_2019}]. This means that making a neural network sparse does not have to degrade performance. 

The sparsity can be introduced to a neural network before, during, or after training. The sparsification in each case is conducted similarly: each parameter is assigned a score based on some metric, and a procedure determines based on the scores whether to keep or remove certain parameters. When sparsifying before or after training, this is only done once. For example, the Optimal Brain Surgeon procedure computes a ``saliency" score based on the Hessian of a fully trained neural network [\cite{hassibi_second_1992}]. If the saliency is much smaller than the loss, then the parameter is removed. When sparsifying during training, this is typically done several times. For example, sparse identification of nonlinear dynamics (SINDy) prunes all neurons with magnitude below 0.1 every 500 epochs [\cite{champion_data-driven_2019}]. Sparsifying before training allows for significant computational savings over the training, but leads to less competitive performance [\cite{gadhikar_why_2023}]. Sparsifying during training is a compromise between computational savings and competitive performance. Sparsifying after training yields no computational savings over training, but can be applied to already trained networks with little impact on performance. 

In reduced-order modelling, sparsity for neural networks mostly comes in two flavours. The first is before training by specifying that the neural network should have a particular structural sparsity. The main example of this is convolutional neural networks, but examples also includes mesh-informed and graph neural networks [\cite{franco_mesh-informed_2023,pichi_graph_2024, carlberg_galerkin_2017}]. In these networks, the weight matrices are no longer dense but instead have many zeroes. The other flavour is through iterative magnitude pruning, through methods like the aforementioned SINDy. Methods from both flavours heavily reduce the effective number of operations required to evaluate the trained networks. However, neither satisfyingly addres what a good choice for the latent dimensionality is, with a common strategy for picking the latent dimensionality being brute force on a what-works-best basis.

A recent work by Bungert et al. showed that the linearized Bregman iterations provide an alternative for creating sparsity in an autoencoder. They abbreviated linearized Bregman iterations to LinBreg, and created an Adam version of it termed AdaBreg. They used both LinBreg and AdaBreg as a replacement for Adam in training tasks with an autoencoder for denoising and deblurring and with a DenseNet for classification [\cite{bungert_neural_2021}]. In their experiments, the autoencoder was given a large architecture with a significant overestimate of the latent dimensionality. Even though the authors did not specify a bottleneck in the architecture, one appeared during training. Whilst this was shown to happen in a computer vision task and not in a reduced-order modelling task, it suggests that these optimizers can help address the question of what a good choice for the latent dimensionality is.

\subsection{Our contribution}
We propose a novel multistep algorithm for creating a sparse autoencoders useful for dimensionality reduction of PDE solution data, with which finding a good latent dimensionality becomes significantly cheaper. The training of the autoencoders is based on the linearized Bregman iterations. In the algorithm, we use the $\norm{\cdot}_{1,2}$ matrix norm and the $\norm{\cdot}_{*}$ nuclear norm to induce maximal sparsity. After training, we perform a latent version of POD to compress the latent space. Subsequently, we clean the network by propagating the biases corresponding to the induced zero rows in weight matrices. The former provides a decrease in latent space dimensionality without meaningfully increasing the loss, whereas the latter converts the sparsity into an effective reduction in the number of operations required to evaluate the autoencoder. The full algorithm is outlined in \cref{alg:sparse-bregman}, and discussed in more detail in \cref{sec:methodology} with background on Bregman iterations in \cref{sec:bregman}.

We applied our method to the solution snapshots of three examples: a 1D diffusion, a 1D advection and a 2D reaction-diffusion equation. The specifics of these examples are detailed in \cref{sec:numerics}. For each of these, we compare the best autoencoders found with our algorithm with Adam and SGD autoencoders. Our autoencoders achieve a similar loss to the reference Adam and SGD autoencoders yet with 30\% fewer parameters and a latent space 60\% smaller. A summary of the results can be found in \cref{tab:summary}.

\subsection{Notation}
In this work, the bold variables are matrices and vectors. The spatial coordinate is $x$ and the vector of nodal points $\vb{x}$. The full-order (FOM) solution is denoted by $\vb{u}_t$, where the lowercase $t$ is the time-index. Derivatives of functions are denoted using the $\partial$ symbol, except that time-derivatives of the latent and full-order state vectors are denoted with a dot above it. The tensor product of two vector spaces $V$ and $W$ is the vector space $V\otimes W$ with a bilinear map $\otimes: V\times W\to V\otimes W$ such that if $B_V$ is a basis for $V$ and $B_W$ is a basis for $W$, then the set $\Set{ v\otimes w \given v\in B_V, w\in B_W}$ is a basis for $V\otimes W$. When we write $\bigotimes_{\ell=1}^L V^\ell$ for a set of vector spaces $\Set{V^\ell}_{\ell=1}^L$, we mean $V^1\otimes\dots\otimes V^L$. $\norm{\vb{W}}_{p,q}:=\bigg(\sum_{j=1}^n\bigg(\sum_{i=1}^m\abs{W_{ij}}^q\bigg)^{p/q}\bigg)^{1/p}$ is the elementwise $\ell^{p,q}$-norm of a matrix $\vb{W}\in \R^{m\times n}$, and $\norm{\vb{W}}_{\ast}:= \sum_{i=1}^{\min\Set{m,n}}\sigma_i(\vb{W})$ with $\sigma_i(\vb{W})$ the $i^{\text{th}}$ singular value of $\vb{W}$ is the nuclear norm of $\vb{W}\in \R^{m,n}$. For a function $f: \X \to \Y$, the Lipschitz constant of $f$ is given by
\begin{equation}
    \lip(f) := \sup_{x,x'\in \X}\frac{\norm{f(x)-f(x')}}{\norm{x-x'}}.
\end{equation}

\section{Bregman iterations}\label{sec:bregman}
\begin{figure}
    \centering
    \begin{subfigure}{0.15\textwidth}
        \includegraphics[width=\textwidth]{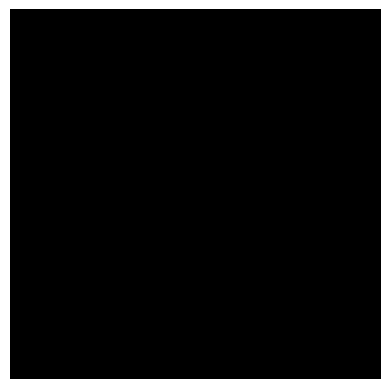}
        \caption{Iteration 0}
    \end{subfigure}    
    \begin{subfigure}{0.15\textwidth}
        \includegraphics[width=\textwidth]{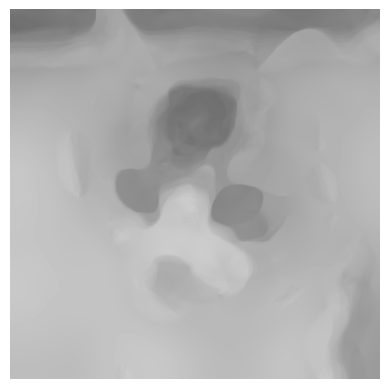}
        \caption{Iteration 1}
    \end{subfigure}    
    \begin{subfigure}{0.15\textwidth}
        \includegraphics[width=\textwidth]{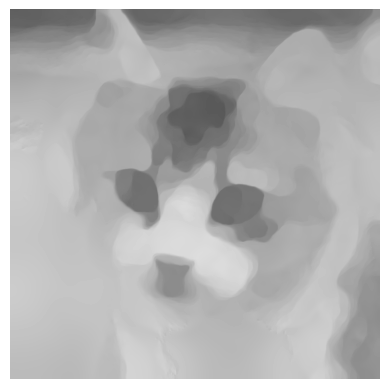}
        \caption{Iteration 2}
    \end{subfigure}    
    \begin{subfigure}{0.15\textwidth}
        \includegraphics[width=\textwidth]{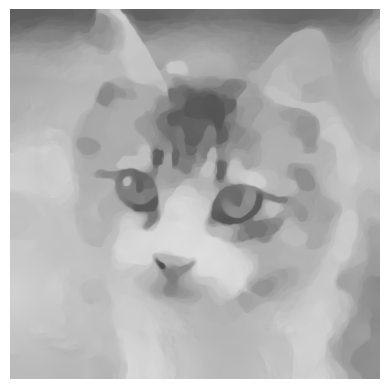}
        \caption{Iteration 4}
    \end{subfigure}    
    \begin{subfigure}{0.15\textwidth}
        \includegraphics[width=\textwidth]{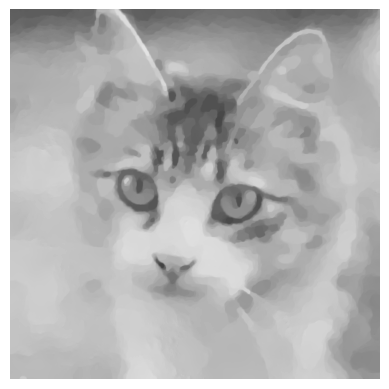}
        \caption{Iteration 8}
    \end{subfigure}    
    \begin{subfigure}{0.15\textwidth}
        \includegraphics[width=\textwidth]{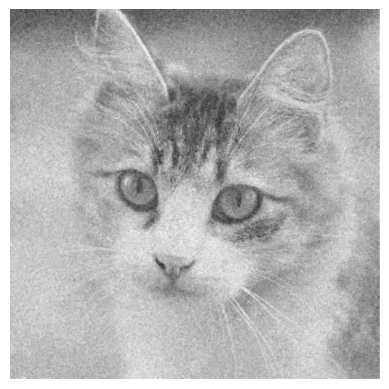}
        \caption{Noisy image}
    \end{subfigure}    
    \caption{Bregman iterations with TV regularization applied to a noisy cat image. The TV regularization is computed using the \textit{PyProximal} [\cite{ravasi_pyproximal_2024}] implementation of the algorithm from \cite{beck_fast_2009} with $\eta=0.5$, $300$ sub-iterations, and a relative tolerance of $\num{1e-8}$.}
    \label{fig:cat}
\end{figure}

The Bregman iterations are the optimization method our multistep algorithm relies on. With the right choice of regularizer it is an example of an algorithm that sparsifies before training. It is developed by Osher et al. based on an algorithm by Bregman [\cite{osher_iterative_2005}]. It solves the bilevel minimization problem
\begin{subequations}\label{eq:bilevel}
\begin{align}
    u^\dagger \in &\argmin_{u\in\mathcal{H}}\RC(u) \\
    \text{s.t. }&u\in\argmin_{\Bar{u}\in\mathcal{H}}\mathcal{L}(\Bar{u}),
\end{align}
\end{subequations}
where $\mathcal{L}$ is a loss function, $\RC$ a regularizer, and $\mathcal{H}$ a Hilbert space [\cite{burger_error_2007}]. The iterations are given by
\begin{equation}\label{eq:bregman_iterations}
    \begin{aligned}
        u_k &= \argmin_{u\in\U} D^{p_{k-1}}_\RC(u,u_{k-1})+\eta\mathcal{L}(u) & u_0 = 0\\
        p_k &= p_{k-1} - \eta\grad\mathcal{L}(u_k) & p_0 = 0, p_k\in\partial \RC(u_k)
    \end{aligned}
\end{equation}
in which the design parameter $\lambda>0$,
\begin{equation}
    \partial\RC(u) =\Set[\bigg]{ p \given \RC(v) \geq \RC(u) + \braket{p}{v-u} \quad \forall v \in \mathcal{H}}
\end{equation}
is the subgradient of $\RC$ at the point $u\in \mathcal{H}$, and 
\begin{equation}
    D^{p}_\RC(u,v) := \RC(u)-\RC(v)-\braket{p}{u-v}
\end{equation}
the Bregman divergence from $u\in\mathcal{H}$ to $v\in\mathcal{H}$ given $p\in\partial \RC(u)$. Note that Bregman divergences are different from metrics, since they do not satisfy the triangle inequality and are not symmetric. The algorithm enjoyed great success in fields like imaging and compressed sensing  [\cite{yin_bregman_2008,bachmayr2009iterative,cai2009linearized,cai2009convergence,yin2010analysis,burger_error_2007,burger_adaptive_2012,benning2018modern,brune_primal_2011,benning_higher-order_2013}]. A simple example shown in \cref{fig:cat}.

Unfortunately, Bregman iterations cannot be applied directly to neural networks, since the loss function with respect to the parameters is in general non-convex. One algorithm that can deal with non-convexity is the linearized Bregman iterations. To go from the Bregman iterations to the linearized Bregman iterations, the regularizer $\RC$ is made strictly convex by adding elastic net regularization $\frac{1}{2}\norm{\cdot}_2^2$ and the loss $\mathcal{L}$ is replaced with its first order Taylor expansion around the previous iterate, i.e. when $\RC$ is replaced with
\begin{equation}
    \RC(\theta) + \frac{1}{2\delta}\norm{\theta}_2^2,
\end{equation}
and $\mathcal{L}$ with
\begin{equation}
    \mathcal{L}(\theta_k) + \braket{\grad \mathcal{L}(\theta_k)}{\theta - \theta_k}.
\end{equation}
for some elastic-net regularization parameter $\delta>0$. Here, $\theta=\Set{\W^L,\b^L,\hdots, \W^1, \b^1}$ are the parameters of the neural network, where $\W^\ell$ and $\b^\ell$ are the weights and biases of layer $\ell$, respectively. After some algebraic manipulations, we get the linearized Bregman iterations
\begin{subequations}\label{eq:linearized_bregman}
\begin{align}
    v_{k+1} &= v_k - \eta\grad \mathcal{L}(\theta_k) \\
    \theta_{k+1} &= \prox_{\delta \RC}(\delta v_{k+1}),
\end{align}    
\end{subequations}
where 
\begin{equation}
    \prox_{\delta \RC}(f) := \argmin_{g\in \mathcal{H}}\frac{1}{2}\norm{f-g}_{\mathcal{H}}^2 + \delta\RC(g) = \argmin_{g\in \mathcal{H}}\frac{1}{2\delta}\norm{f-g}_{\mathcal{H}}^2 + \RC(g)
\end{equation}
is the proximal operator for $\delta \RC$ over a Hilbert space $\mathcal{H}$ with a chosen elastic-net regularization parameter $\delta>0$ [\cite{bungert_bregman_2021}]. This proximal operator is a generalization of a projection and appears in many optimisation algorithms [\cite{parikh_proximal_2014}]. Although it seems that in every iteration of \cref{eq:linearized_bregman} we need to solve a minimization problem because of the proximal operator, it turns out that the proximal operator has a closed form expression for many regularizers $\RC$. Linearized Bregman iterations are a generalization of mirror descent, which in turn is a generalization of gradient descent [\cite{beck_mirror_2003}]. They are also closely related to proximal gradient descent with the main difference being that linearized Bregman apply the proximal operator to the accumulated gradients instead of only to the current one [\cite{benning_choose_2021}]. Note also that this scheme has 3 parameters: the elastic-net regularization parameter $\delta$, the learning rate $\eta$ and the regularization parameter $\lambda$. Since the factors inside a minimization problem can be multiplied by a positive constant $c>0$ without changing the outcome, one of these parameters is redundant. In particular, $(\delta, \eta, \lambda)$ will yield the same iterations as $(c\delta, \eta/c, \lambda/c)$ for any $c>0$. Hence, we can choose $\delta=1$ and vary $\eta$ and $\lambda$ without loss of generality.

Linearized Bregman iterations have been used in \cite{bungert_neural_2021} with autoencoders for denoising and deblurring and with a DenseNet for classification. In \cite{bungert_bregman_2021}, the Adam and the momentum version of linearized Bregman iterations have been developed. They abbreviated Linearized Bregman iterations to LinBreg and called the Adam version AdaBreg. These works show that using these algorithms with the sparsity-inducing regularizer
\begin{equation}
    \RC(\theta) := \lambda\sum_{\ell=1}^L\sqrt{d^\ell}\norm{\W^\ell}_{1,2}
\end{equation}
gives networks that perform similarly to Adam yet with significantly fewer parameters and the appearance of a bottleneck. During this training, the density, latent size and accuracy of the networks increased over the epochs. This is similar to how adding POD modes increases the accuracy of reconstruction and number of parameters used in this projection. We adapt the full training procedure for our reduced-order modelling tasks in the following section.

\section{Methodology}\label{sec:methodology}
In this work, we consider the semi-discretized form of a partial differential equation given by the high-dimensional dynamical system
\begin{equation}\label{eq:full_order_dynamics}
\begin{aligned}
    \Dot{\vb{u}}_t &= f(\vb{u}_t,\mu) \\
    \vb{u}_0 &= g(\mu)
\end{aligned}
\end{equation}
with time $t\in[0,T]$, parameters $\mu$ relevant for this problem, and (non)linear functions $f$ and $g$. We seek to find an appropriate latent space for the solutions to \cref{eq:full_order_dynamics}. This latent dimensionality should be 
\begin{itemize}
    \item small but 
    \item large enough to capture the relevant information from the full-order solution $\vb{u}_t$ and 
    \item the maps from full-order space to latent space and vice versa should be easy to compute.
\end{itemize}
For the maps we use an autoencoder $\phiAE(\vb{x})=(\phiDec\circ \phiEnc)(\vb{x})$ where both the decoder $\phiDec$ and the encoder $\phiEnc$ are simple feedforward dense neural networks. These are  parameterised by $\theta=\Set{\W^L,\b^L,\hdots, \W^1, \b^1}$, where $\W^\ell$ and $\b^\ell$ are the weights and biases of layer $\ell$, respectively. We use the Linearized Bregman iterations to find the parameters that achieve the aforementioned goals. In particular, we choose 
\begin{subequations}
\begin{align}
    \mathcal{H} &= \bigotimes_{\ell=1}^L F^\ell \otimes \ell^2(d^\ell) \\
    \mathcal{L}(\bm{\theta}) &= \sum_{\vb{u}\in\mathcal{X}}\norm{\vb{u}-\phiAE(\vb{u})}^2_2, \\
    \RC(\theta) &= \lambda\bigg(\sum_{\substack{\ell=1\\ \ell\neq L_{\text{enc}}}}^L\sqrt{d^\ell}\norm{\W^\ell}_{1,2} + \norm{\W^{L_{\text{enc}}}}_\ast\bigg).
\end{align}
\end{subequations}
as the Hilbert space, loss function, and regularizer, respectively, in which $d^\ell$ is the number of rows of $\W^\ell$ or equivalently the length of $\b^\ell$ with $d^0$ the full-order dimensionality, $F^\ell$ is the Hilbert space of $d^\ell\times d^{\ell-1}$ matrices with the Frobenius norm, $\otimes$ is the tensor-product between vector spaces, $L_{\text{enc}}$ refers to the last layer of the encoder, $\norm{\cdot}_{1,2}$ is the matrix $\ell^{1,2}$ norm, $\norm{\cdot}_*$ the nuclear norm, and $\mathcal{X}$ is a training dataset consisting of solutions to \cref{eq:full_order_dynamics} at several time instances for parameter values. The snapshot matrix $\vb{X}$ is the matrix with as columns $\vb{u}\in\mathcal{X}$. The chosen loss function ensures the training algorithm captures the relevant information (second goal), whereas the elementwise norms in the regularizer ensure the encoder and decoder are sparse (third goal) and the nuclear norm ensures the latent dimensionality is small (first goal). By making the Hilbert space a direct sum space and choosing the regularizer to be a sum of norms of weights or biases, we can use the separation property for proximal operators to control the impact of the proximal on each weight matrix independently. 

\textbf{Separation property:} Let the Hilbert space $\mathcal{H}=\otimes_{i=1}^N\mathcal{H}_i$ be a direct sum of $N$ Hilbert spaces $\mathcal{H}_i$. If $f\in C(\mathcal{H})$ with $f(x) = \sum_{i=1}^Nf_i(x)$ and $f_i\in C(\mathcal{H}_i)$ convex, then 
\begin{equation}
    \prox_f(A) = \mqty( \prox_{f_1}(A_1) & \hdots & \prox_{f_N}(A_N))
\end{equation}
for all $A\in \mathcal{H}$.

In particular, we can apply the proximal operator with $\norm{\cdot}_{1,2}$ to the weight matrices $\W^\ell$  with $\ell\neq L_{\text{enc}}$, the proximal with $\norm{\cdot}_*$ to the weight matrix $\W^{L_{\text{enc}}}$ and the proximal with no regularizer to biases $\b^\ell$. This means the linearized Bregman iterations for biases $\b^\ell$ reduces to stochastic gradient descent, since
\begin{equation}
    \prox_{0}(\b^\ell) = \argmin_{\b \in \ell^2(d^\ell)}\frac{1}{2}\norm{\b^\ell-\b}^2_{2} = \b^\ell
\end{equation}
for all $\ell$.

\begin{figure}
    \centering
    \begin{subfigure}{0.3\textwidth}\label{fig:methodology_a}
        \includegraphics[width=\textwidth]{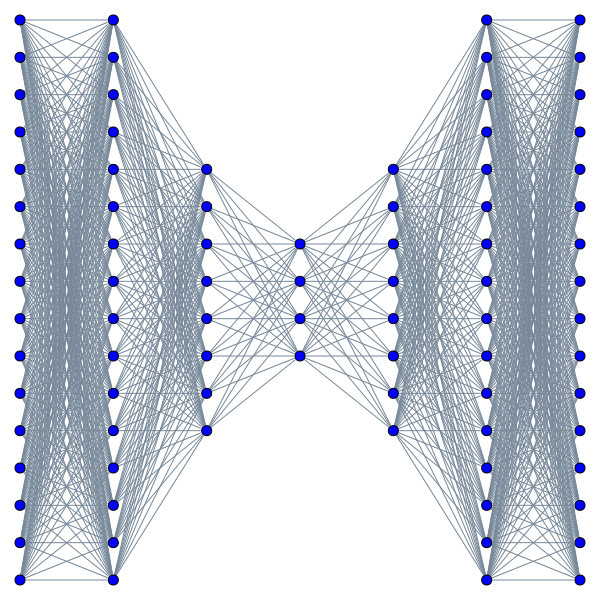}
        \caption{Initialized network}
    \end{subfigure}
    \hfill
    \begin{subfigure}{0.3\textwidth}
        \includegraphics[width=\textwidth]{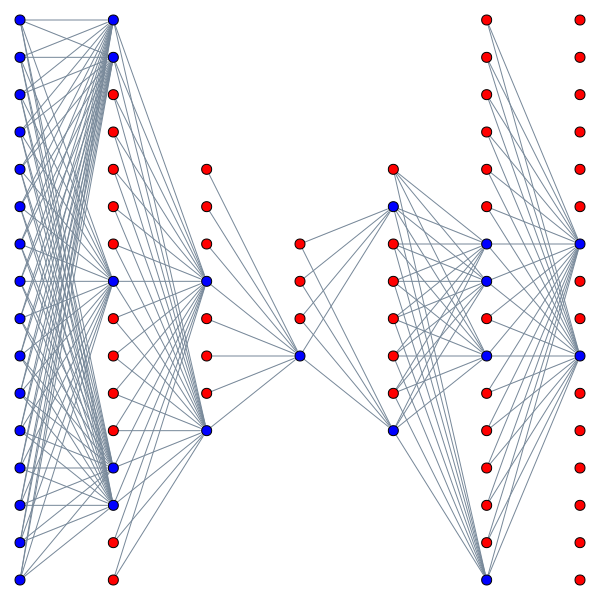}
        \caption{Sparsified network}
    \end{subfigure}
    \hfill
    \begin{subfigure}{0.3\textwidth}
        \includegraphics[width=\textwidth]{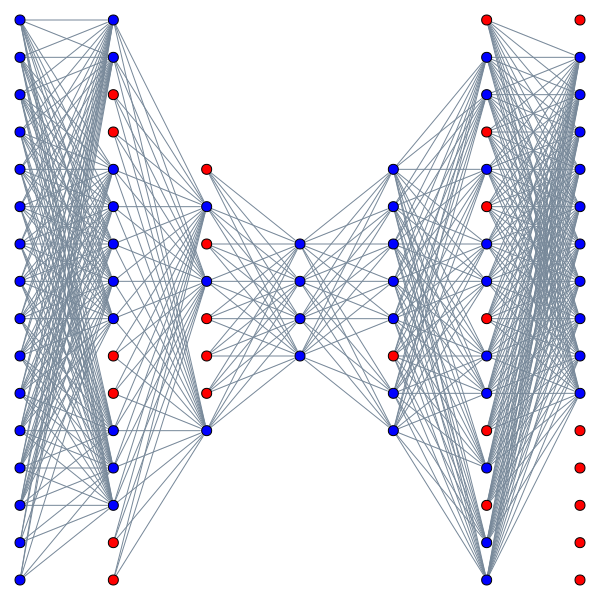}
        \caption{Trained network}
    \end{subfigure}
    \\ 
    \begin{subfigure}[t]{0.3\textwidth}
        \includegraphics[width=\textwidth]{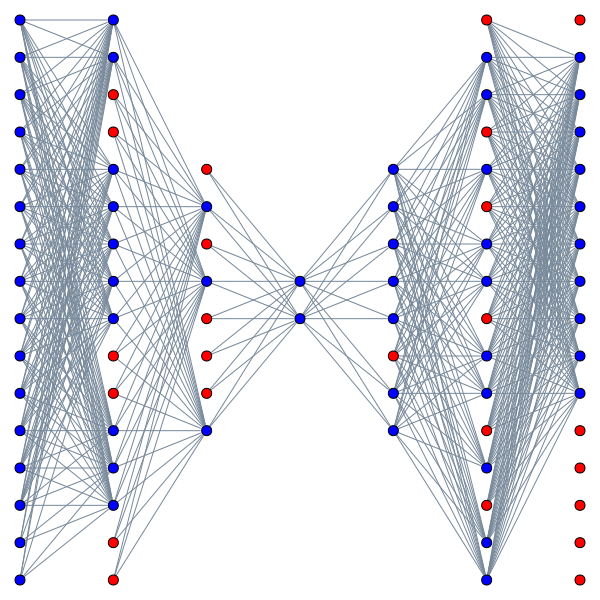}
        \caption{Network after latent truncated SVD}
    \end{subfigure}%
    ~
    \begin{subfigure}[t]{0.3\textwidth}
        \includegraphics[width=\textwidth]{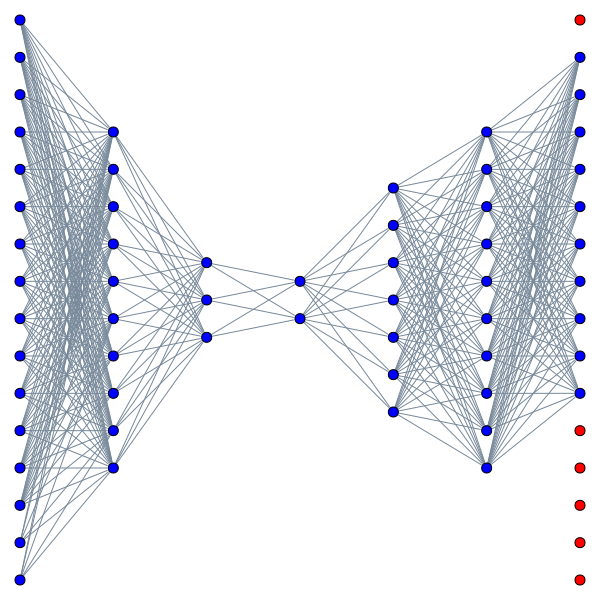}
        \caption{Network after propagation}
    \end{subfigure}
    \caption{Example of \cref{alg:sparse-bregman} applied to the small autoencoder. Blue nodes represent active neurons (output still depends on input), red inactive neurons (output does not depend on input), and lines between two neurons indicate the weight between them is non-zero. In (a), we see the autoencoder after lines 1 and 2 of \cref{alg:sparse-bregman} have been executed. In (b), we see the autoencoder after lines 3 to 9 of \cref{alg:sparse-bregman} have been executed. A lot of inactive neurons have been created, but we note by construction each of these is still outputting a non-zero value. This is due to the biases $b^\ell$ being strictly positive. In (c), we see the autoencoder after training. Roughly half of the inactive neurons have been become active. In (d), we see the network after having applied the latent SVD post-processing described in lines 20 to 25 of \cref{alg:sparse-bregman}. The difference between (c) and (d) is that half the neurons in the bottleneck have been removed. In (e), we see the network after having applied the bias propagation post-processing described in lines 26 and 27. Visually, this looks like removing inactive neurons from all layers except from the output. Note that in this particular example, there are still inactive neurons in the output layer. This is caused by the data all having the same  values for the indices corresponding to these neurons. This implies that the training algorithm is also able to detect this kind of pattern in the data.}
    \label{fig:methodology}
\end{figure}

In \cref{sec:initialization} we discuss how to initialize the neural network and dual variable $v^{(0)}$ during training runs, in \cref{sec:truncated} we discuss the reasoning for the nuclear norm in more detail, and in \cref{sec:propagation} we discuss the reasoning for the weighted $\ell^{1,2}$ norm in more detail. Pseudocode for the entire algorithm when using LinBreg is shown in \cref{alg:sparse-bregman}. When AdaBreg instead of LinBreg, lines 14 and 15, which refer to the LinBreg iterations, will be replaced with corresponding AdaBreg iterations from \cite{bungert_bregman_2021}. An example of \cref{alg:sparse-bregman} applied to a small autoencoder is visualized in \cref{fig:methodology}.

\begin{algorithm}[hbt!]
\caption{Sparse Bregman for LinBreg}\label{alg:sparse-bregman}
\KwData{learning rate $\eta$, regularization constant $\lambda$, initial density level $p$, initialized autoencoder $\bm{\phi}_\theta$, truncation level $\epsilon$}
\KwResult{Parameters for the sparse autoencoder $\theta$}

\BlankLine  \Comment*[r]{initialization}
$\W^\ell \sim \mathcal{U}(-\sqrt{6}/d^{\ell-1},\sqrt{6}/d^{\ell-1})$\;
$\b^\ell \sim \mathcal{U}(1,1/d^{\ell-1})$

\BlankLine  \Comment*[r]{Sparsify the initialization}
$\vb{U}, \vb{S}, \vb{V}^\intercal\gets $ singular values decomposition of $\W^{L_{\text{enc}}}$\;
$\W^{L_{\text{enc}}} \gets \vb{U} \diag(1, 0, \hdots, 0)\vb{V}^\intercal$
\For{$\W^\ell \in \theta$ with $\ell\neq L_{\text{enc}}$}{
    $S\gets \lceil d^\ell/p \rceil$ unique random integers $i \sim U\Set{1,d^\ell}$\;
    \For{$i\in S$}{
        replace the $i^{\text{th}}$ row of $\W^\ell$ with a row of zeroes
    }
}
\BlankLine \Comment*[r]{Train}
\For{each epoch}{
    \For{each batch in epoch}{
        Evaluate $\mathcal{L}(\theta_k)$ over batch\;
        Backpropagate to compute $\grad\mathcal{L}(\theta_k)$\;
        $v_{k+1}=v_k-\eta\grad\mathcal{L}(\theta_k)$\;
        $\theta_{k+1} = \prox_{\RC}(v_{k+1})$\;
        $k \gets k+1$\;
    }
}
\BlankLine \Comment*[r]{Latent SVD}
$\vb{U}, \vb{S}, \vb{V}^\intercal\gets $ singular values decomposition of $\W^{L_{\text{enc}}}$\;
$\vb{U}, \vb{S}, \vb{V}^\intercal\gets $ truncation of $\vb{U}, \vb{S}, \vb{V}^\intercal$ at error level $\epsilon$\;
$\W^{L_{\text{enc}}} \gets \vb{SV}^\intercal $\;
$\W^{L_{\text{enc}}+1} \gets \W^{L_{\text{enc}}+1}\vb{U}$\;
$\b^{L_{\text{enc}}} \gets 0$\;
$\b^{L_{\text{enc}}+1} \gets \W^{L_{\text{enc}}+1} \b^{L_{\text{enc}}}+\b^{L_{\text{enc}}+1}$\;
\BlankLine \Comment*[r]{Propagation}
propagate bias on $\bm{\phi}_\theta$ using the algorithm from the \textit{Simplify} package\;
restore $\W^L$
\BlankLine 
\end{algorithm}

\subsection{Sparse initialization}\label{sec:initialization}
During the training with LinBreg or AdaBreg, the network will become increasingly dense. Hence, we want to initialize it sparse.  On one hand, it is unlikely that the network becomes sparser than the initialized density during training. On the other, initializing the network too sparse implies limited information is passed through the early stages of training. Hence, we will initialize the networks with 20\% of the rows of weights $\W^\ell$ non-zero for all $\ell\neq L_{\text{enc}}$, and the biases dense. We also initialise $\W^{L_{\text{enc}}}$ spectrally sparse by truncating all singular values except for the largest one. 

More precisely, to achieve the desired sparse initialization, we first initialize a dense network with Kaiming uniform initialization. The network uses ReLU as activation function. Hence, this means that for the weights we take $\W^\ell\sim \mathcal{U}(-\sqrt{6}/d^{\ell-1},\sqrt{6}/d^{\ell-1})$ and for the biases $\b^\ell\sim\mathcal{U}(0,1/d^{\ell-1})$. Here, the uniform interval for the biases is taken over the positive part of the interval instead of the interval symmetric around zero to ensure that regardless of the sparsification of $\W^\ell$ the neurons will still give a non-zero output. Subsequently, for the weight matrices $\W^\ell$ with $\ell\neq L_{\text{enc}}$ we compute the number of rows for each matrix that need to be zeroed to achieve 20\% sparsity using $\lceil d^\ell / (1-0.2) \rceil$. We sample that many unique integers uniformly from 1 to $d^\ell$. For each of the sampled integers, we make the corresponding row of $\W^\ell$ into a zero row. For the weight matrix $\W^{L_{\text{enc}}}$, we compute its SVD $\W^{L_{\text{enc}}}=\vb{USV}^\intercal$ and replace it with $\vb{U}\diag(1, 0, \hdots, 0)\vb{V}^\intercal$. We do not sparsify the biases $\b^\ell$.

\subsection{Latent truncated SVD}\label{sec:truncated}
The nuclear norm $\norm{\W^{L_{\text{enc}}}}_*$ in the regularizer $\RC$ enables an SVD-based post-processing of the network. This post-processing step works best when the weight matrix $\W^{L_{\text{enc}}}$ has a few large singular values and the remaining ones are small, analogously to the POD. The nuclear norm $\norm{\W^{L_{\text{enc}}}}_*$ aides in achieving this. 

For an arbitrary matrix $\A\in F^{L_{\text{enc}}}$, it has the associated proximal term
\begin{equation}
    \prox_{\lambda \norm{\cdot}_*}(\A) = \argmin_{\B\in F^{L_{\text{enc}}}}\frac{1}{2\lambda}\norm{\A-\B}^2_{F^{L_{\text{enc}}}}+\norm{\B}_* = \vb{U} \diag([s_1,s_2,...]) \vb{V}^\intercal
\end{equation}
with $\A=\vb{USV}^\intercal$ the singular value decomposition of $\A$ and 
\begin{equation}
    s_i = \begin{cases}
        0 & S_{ii} \leq \lambda, \\
        S_{ii} - \lambda & S_{ii} \geq \lambda.
    \end{cases}
\end{equation}
This shifts and truncates the singular values of $A$, encouraging a few larger singular values and the remaining singular values to be small or zero. Hence, we get the desired singular values for $\W^{L_{\text{enc}}}$ for which the method works best when adding $\norm{\W^{L_{\text{enc}}}}_*$ to the regularizer.

The post-processing works by using the fact that we have no activation function between $\W^{L_{\text{enc}}}$ and $\W^{L_{\text{enc}}+1}$. Since there is no activation function, we can replace the matrices $\W^{L_{\text{enc}}}$ and $\W^{L_{\text{enc}}+1}$ by any matrices $\A$ and $\B$ that satisfy
\begin{equation}
    \W^{L_{\text{enc}}+1}\W^{L_{\text{enc}}} = \A\B
\end{equation}
without affecting the full autoencoder $\bm{\phi}_\theta$. Note that it does impact the latent space. In particular, we can take the singular value decomposition of $\W^{L_{\text{enc}}}$ so that
\begin{equation}
    \W^{L_{\text{enc}}+1}\W^{L_{\text{enc}}} = \W^{L_{\text{enc}}+1}(\vb{USV}^\intercal) =  (\W^{L_{\text{enc}}+1}\vb{U})(\vb{SV}^\intercal).
\end{equation}
This allows us to create a new autoencoder $\bm{\phi}_{\theta_{\text{post}}}$ from the old autoencoder $\bm{\phi}_{\theta}$ using the mapping
\begin{subequations}    
\begin{align}
    \W^{L_{\text{enc}}} &\mapsto \vb{SV}^\intercal, \\
    \W^{L_{\text{enc}}+1} &\mapsto \W^{L_{\text{enc}}+1}\vb{U}, \\
    \b^{L_{\text{enc}}} &\mapsto 0, \\
    \b^{L_{\text{enc}}+1} &\mapsto \W^{L_{\text{enc}}+1} \b^{L_{\text{enc}}}+\b^{L_{\text{enc}}+1}.
\end{align}
\end{subequations}
Both $\bm{\phi}_{\theta_{\text{post}}}$ and $\bm{\phi}_{\theta}$ have the same reconstruction error. When $\W^{L_{\text{enc}}}$ is replaced with its SVD truncated version, the size of the latent space would not change. By moving the $U$ matrix to the decoder, there is a change upon replacing $\vb{U}$, $\vb{S}$ and $\vb{V}$ with their truncated versions $\vb{U}_\epsilon$, $\vb{S}_\epsilon$ and $\vb{V}_\epsilon$. The size of the truncated latent space is determined by the number of singular values kept in $\vb{S}_\epsilon$. The decoder has to be taken into account when determining the number of singular values to keep. It is Lipschitz continuous when the activation function $\sigma$ used is Lipschitz. Hence, truncating at an error level $\epsilon$ will introduce an $\ell^2$-error of at most $\lip(\bm{\phi}_{\theta_{\text{post}},\text{dec}})\epsilon$. Choosing $\epsilon$ small enough (e.g. 1\% of the reconstruction error) will means that other errors dominate. Hence, the latent dimension can be reduced effectively without worsening the result.

\subsection{Bias propagation}\label{sec:propagation}
The term $\norm{\W^\ell}_{1,2}$ in the regularizer enables another kind of post-processing. During this post-processing, useless zero rows are removed from the weight matrices. The terms $\norm{\W^\ell}_{1,2}$ encourage the creation of these zero rows.

For an arbitrary matrix $\A\in F^\ell$, the term $\norm{\cdot}_{1,2}$ has the associated proximal term
\begin{equation}
    \prox_{\lambda\norm{\cdot}_{1,2}}(\A)_{ij} = \bigg(\argmin_{\B\in F^\ell}\frac{1}{2\lambda}\norm{\A-\B}^2_{F^\ell}+\norm{\B}_{1,2}\bigg)_{ij} = \begin{cases}
        0 & \norm{\A_{:,j}}_{2} \leq \lambda, \\
        \abs{A_{ij}} - \sign(A_{ij})\lambda & \norm{\A_{:,j}}_{2} \geq \lambda.
    \end{cases} 
\end{equation}
This makes rows of $\W^\ell$ zero when the $\norm{\cdot}_{\ell^2}$ norm of the row falls below a certain threshold. Above the threshold, it shifts the values in the row towards zero. Since pre-activations have the form $\W^\ell x^\ell + \b^\ell$, we have that the $i^{\text{th}}$ output is independent of $x^\ell$ when the $i^{\text{th}}$ row of $\W^\ell$ is zero. This also implies that the $i^{\text{th}}$ column of $\W^{\ell+1}$ is independent of $x^\ell$. Hence, $\W^{\ell}_{i,:}\sigma(\b^\ell_i)$ can be added to the bias of the next layer and the $i^{\text{th}}$ row of $\W^\ell$, the $i^{\text{th}}$ column of $\W^{\ell+1}$, and $\b^\ell_i$ removed from the network without affecting the full autoencoder $\phiAE$. By propagating the biases like this, we effectively reduce the number of parameters in the network. This idea was originally used in \cite{bragagnolo_simplify_2022} for convolutional neural networks and made public using their PyTorch package \textit{Simplify}. Initially, this package did not support the networks used in this paper. We collaborated with them to expand the functionality of their package to provide support for these networks. 

Special care has to be taken with the final layer of the autoencoder. In this layer, we don't want to remove the zero rows. Doing so, will change the output dimensionality of the autoencoder. The Simplify package supports an option that makes sure that this does not happen. Their implementation uses a special layer in that case for the output layer, which uses `torch.scatter` to construct the output of the network. This is inefficient for a network where many of the outputs are dependent on the input, like is the case with autoencoders. In some basic tests, it was shown to be more efficient to turn the final layer back into a simple linear layer.

\section{Numerical experiments}\label{sec:numerics}
In this section, we compare \cref{alg:sparse-bregman} given with LinBreg and its AdaBreg version with SGD and Adam on data produced by a 1D diffusion equation, 1D advection equation and 2D reaction diffusion equation. We motivate the choice of these equations at the start of each subsection, and show how we solved them numerically. The optimizers all need hyperparameters with different orders of magnitude. To determine the best hyperparameters, we use \textit{Weights and biases} (Wandb) hyperparameter optimization \textit{sweeps}. During a sweep, various networks will be trained based on the settings provided to the sweep. In our experiments, these settings include the architecture of the autoencoder, the optimizer to use during training, the number of epochs, batch size and learning rate. For each optimizer and equation combination, we run two sweeps when the optimizer is either SGD or Adam and five sweeps when the optimizer is LinBreg or AdaBreg. We compare the best model from the last sweep for each optimizer on the test data.

During all sweeps but the last, we fix all but the learning rate to fixed values. For the learning rate, we provide a lower and an upper bound. The four sweeps for LinBreg and AdaBreg have a different but fixed regularization constant. The first two runs of each sweep have a learning rate sampled uniformly from the given range. The runs after that in the same sweep use a Gaussian process and acquisition function to determine the next learning rate. First, a Gaussian process using a Matérn 3/2 kernel is fitted to the learning rate and loss tuples. Then, the acquisition function is used to compute the learning rate to use. The acquisition function is based on  \textit{expected improvement}. It assigns a value based on the current lowest loss and the mean and variance of the Gaussian process to 1000 sampled learning rates from the given range. The chosen learning rate is the one with the highest value among these points. 

During the last sweep of each case, the same process in the previous paragraph is repeated but with slight changes. For each optimizer and equation, we determined the best regularization constant and learning rate combination based on the earlier sweeps. In these last sweeps, the regularization constant is fixed to the best regularization constant and the range of the learning rate is set narrowly around the best learning rate. Each sweep yields a best model. For each equation, we compare the best model for each optimizer with each other. 

\subsection{1D diffusion}\label{sec:diffusion}
\begin{figure}
    \centering
    \includegraphics[width=0.48\textwidth]{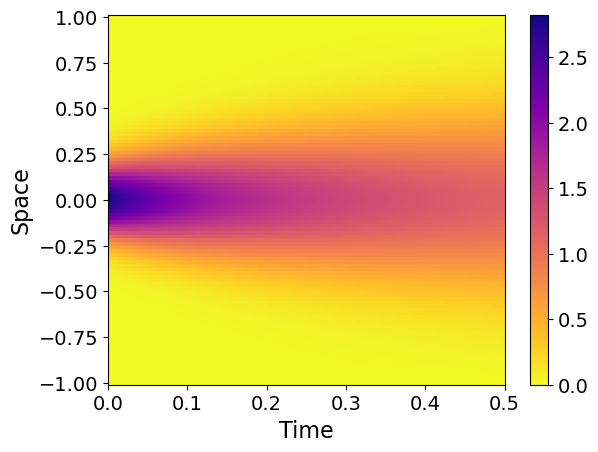}
    \includegraphics[width=0.48\textwidth]{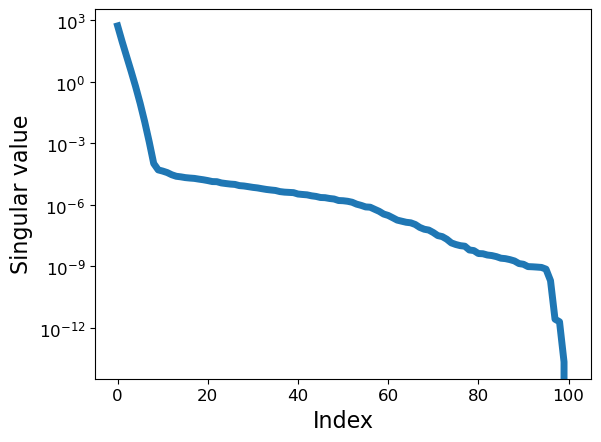}
    \caption{Numerical solution on the left for the diffusion equation with $\mu_{\text{diff}}=0.1$, and  on the right the singular value decay of the snapshot matrix corresponding to the shown numerical solution.}
    \label{fig:diffusion}
\end{figure}
The first example we consider is the 1D diffusion equation, which is given by 
\begin{subequations}\label{eq:diffusion_full}
    \begin{align}
         \partial_t u(x,t) &= \mu_{\text{diff}}\partial_{xx} u(x,t) \\
         u(-1,t) &= u(1,t) = 0 \\
         u(x,0) &= g(x)
    \end{align}
\end{subequations}
on the interval $x\in I=[-1,1]$ for $t\in [0,T]$, where 
\begin{equation}\label{eq:initial_condition_diff}
    g(x) = \frac{1}{\sqrt{0.04\pi}}e^{-\frac{x^2}{0.04}},
\end{equation}
i.e. a Gaussian with mean $0.5$ and variance $0.02$. The exact solution of \eqref{eq:diffusion_full} is given by 
\begin{equation}\label{eq:diffusion_solution}
\begin{aligned}
    u(x,t) &= \sum_{n=1}^\infty \braket{\phi_n}{g}\phi_n(x)\e^{-\pi^2n^2\mu_{\text{diff}}t },
\end{aligned}
\end{equation}
where $\phi_n(x) = \sin(n\pi x)$. The presence of the $n^2$ in the exponential hints that a basis of the first few $\phi_n$ can approximate $u$ already quite well. A full-order solution can be seen in the left figure of \cref{fig:diffusion}. The right figure shows the singular values of the snapshot matrix. The rapid decay of the singular values agrees with that only a few modes are sufficient for accurate projection. Therefore, applying the Bregman framework to this snapshot matrix allows us to provide a comparison in a case that POD works great. 

For the FOM, the interval $I=[-1,1]$ was divided into $N_x=101$ nodal points and $T=1$ taken with $N_t=5001$ time steps. This implies the time step size $\Delta t$ and grid size $\Delta x$ are given by $\Delta t=\num{2e-4}$ and $\Delta x=0.01$, respectively. Every trajectory was subsampled by storing only every $20^{\text{th}}$ time step. For the diffusion parameter, we took $\mu_{\text{diff}}\in \Set{0.1, 0.5, 1}$ for the training set and $\mu_{\text{diff}}\in \Set{0.6}$ for the testing set. To evolve the wave in time, we used a finite difference scheme. In particular, 
\begin{equation}
    \vb{u}_{n+1,i} = (1-2c)\vb{u}_{n,i}+ c(\vb{u}_{n,i-1}+\vb{u}_{n,i+1})
\end{equation}
for all $i\in \Set{0,\hdots,100}$ and $\vb{u}_{n,0}=\vb{u}_{n,100}=0$ for all $n$, where $c=\mu_{\text{diff}}\frac{\Delta t}{(\Delta x)^2}$. The solution with $\mu_{\text{diff}}=0.1$ is shown in \cref{fig:diffusion} together with the singular values of the snapshot matrix.

\subsubsection*{Architecture selection}
When POD is applied to the training set, $r=5$ singular values are required to get to the relative reconstruction error over the training set below $10^{-5}$. Hence, the autoencoders are set to have layers sizes $(101, 50, 25, 5, 25, 50, 101)$. 

For both training and testing, a batch size of $64$ is used. We train for $\num{5000}$ epochs.

\subsubsection*{Hyperparameter optimisation}
To find the best hyperparameters, four experiments have been conducted: one for each optimizer. These require one sweep for SGD as well as Adam and four sweeps for LinBreg as well as AdaBreg.

\begin{figure}[h!]
    \centering
    \begin{minipage}{0.48\textwidth}
        \centerline{\includegraphics[width=\textwidth]{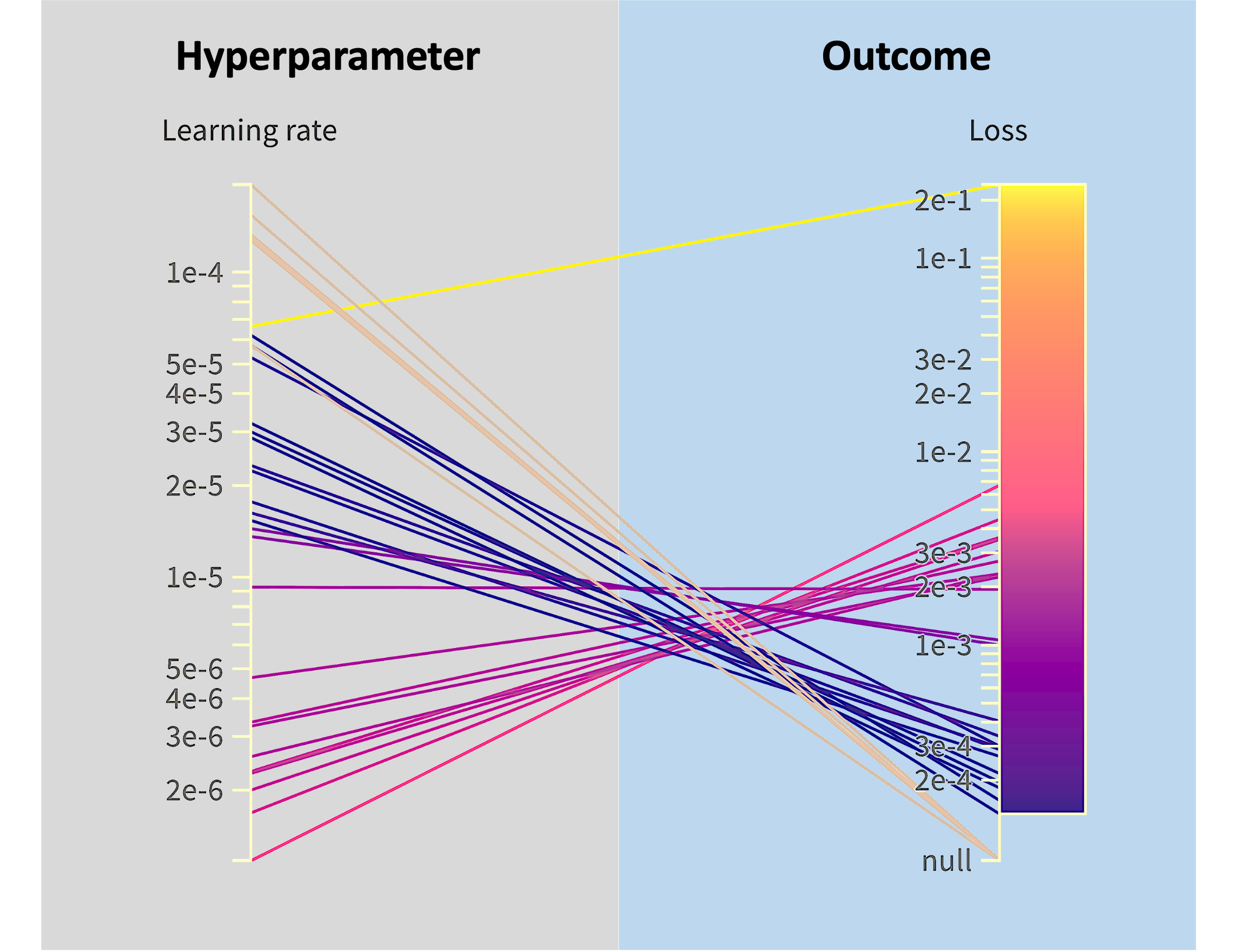}}
        \caption{Wandb sweep for SGD. Runs are started with the learning rate on the left, and the lowest training loss achieved is shown on the right. The lowest loss is achieved for $\eta\approx\num{2e-5}$.}
        \label{fig:diffusion_sgd_sweep}
    \end{minipage}
    \begin{minipage}{0.48\textwidth}
        \centerline{\includegraphics[width=\textwidth]{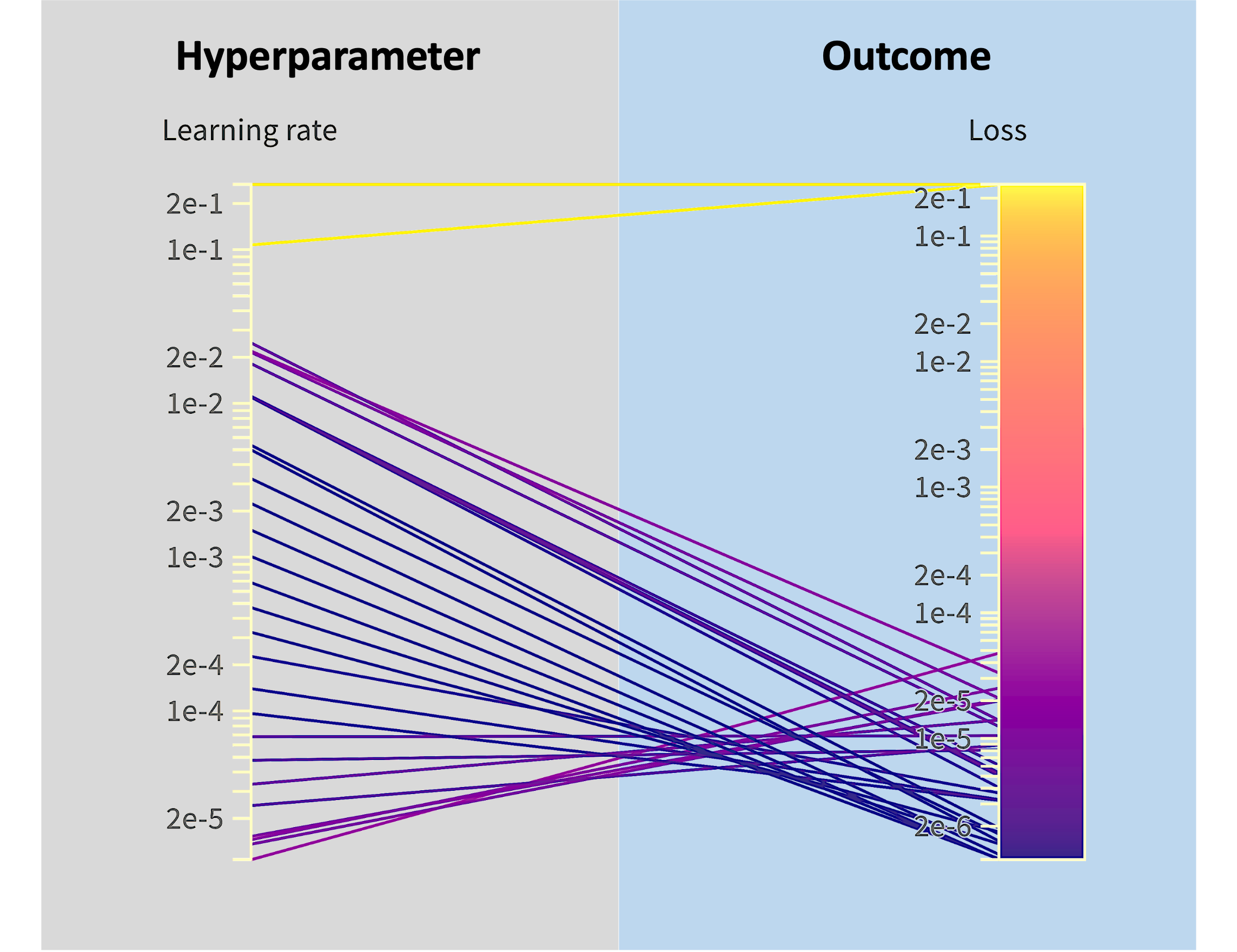}}
        \caption{Wandb sweep for Adam. Runs are started with the learning rate on the left, and the lowest training loss achieved is shown on the right. The lowest loss is achieved for $\eta\approx\num{1.5e-3}$.}
        \label{fig:diffusion_adam_sweep}
    \end{minipage}
\end{figure}

The first optimizer considered is SGD. Here, the prior runs from $\eta\approx \num{e-6}$ to $\eta\approx \num{2e-4}$. The runs are shown in \cref{fig:diffusion_sgd_sweep}. For learning rates below $\eta=\num{5e-6}$, training is too slow due to the learning rate being too small. For learning rates above $\eta=\num{e-4}$, there are a lot of crashed runs. For $\eta\approx \num{5e-5}$, the oscillations are still small. Higher learning rates lead to exploding gradients, so we take $\eta=\num{5e-5}$ as the best learning rate.

The second optimizer considered is Adam. In this case, the prior runs from $\eta\approx \num{2e-5}$ to $\eta\approx \num{2e-1}$. The runs are shown in \cref{fig:diffusion_adam_sweep}. From $\eta=\num{2e-5}$ to $\num{1.5e-3}$ the loss decreasing, after which it starts increasing again. Hence, we take $\eta=\num{1.5e-3}$ as the best learning rate. 

\begin{figure}[h!]
    \centerline{
    \includegraphics[width=1.1\textwidth]{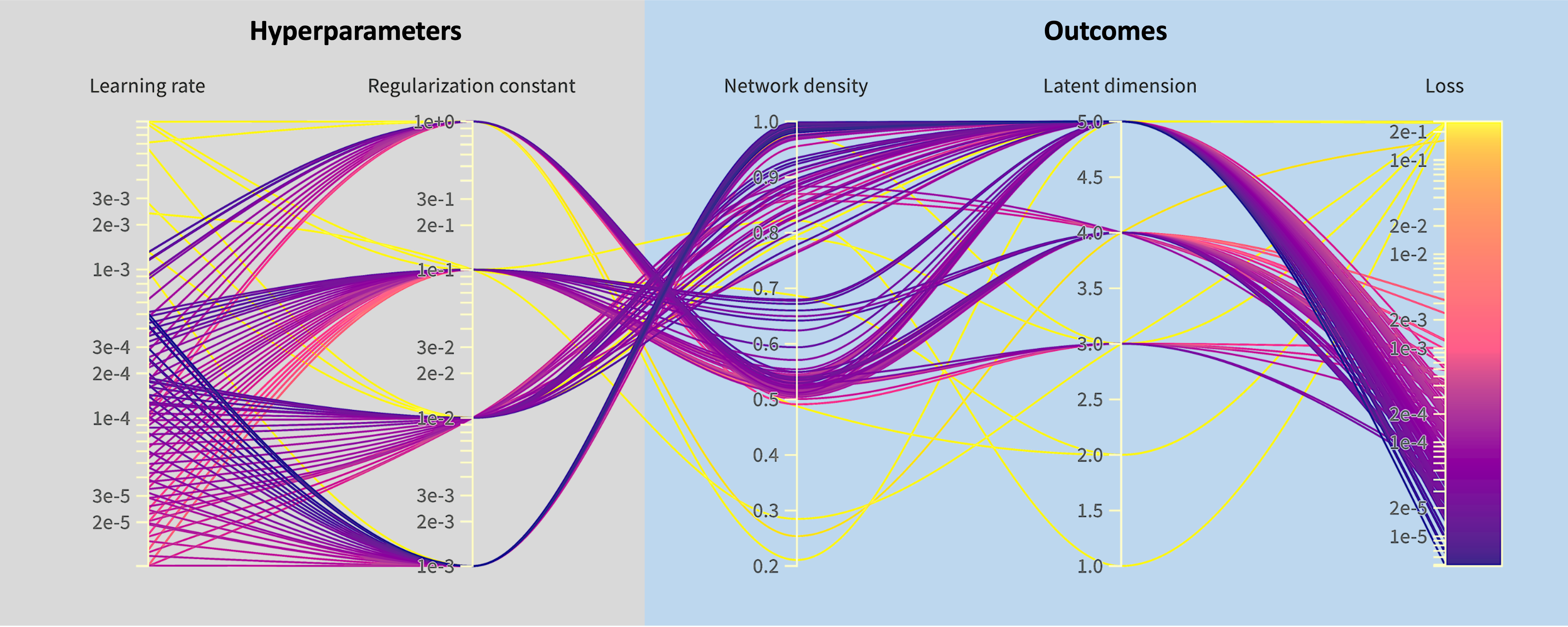}
}   
\caption{Results for Wandb sweeps for the 1D diffusion equation and LinBreg. Runs are started with the hyperparameters on the left, and the outcomes for the network with the lowest training loss achieved are shown on the right.}
    \label{fig:diffusion_linbreg_sweep}
\end{figure}

For LinBreg, four different sweeps are performed. For each of those sweeps, a different regularization constant $\lambda$ was used. These were taken one order of magnitude apart, i.e. sweep $i$ used $\lambda = 10^{i-1}$. The learning rates for these sweeps are taken with bounds $\eta=\num{1e-4}$ till $\eta=\num{1e-2}$, $\eta=\num{1e-5}$ till $\eta=\num{1e-2}$, $\eta=\num{1e-5}$ till $\eta=\num{1e-2}$, and $\eta=\num{1e-5}$ till $\eta=\num{1e-3}$, respectively. The results of these sweeps are visualized in \cref{fig:diffusion_linbreg_sweep}. The lowest training losses are achieved with $\eta=\num{5e-4}$ with $\lambda=0.001$, but for those runs the network is dense and the latent space full. A $20\%$ reduction in latent dimension and a $50\%$ reduction in the number of parameters can be gained at the cost of one order of magnitude in the loss by choosing $\lambda=1$ and $\eta\approx \num{e-3}$. Hence, we take as learning rate $\eta=\num{e-3}$ and as regularization constant $\lambda=1$.

\begin{figure}[h!]
    \centerline{\includegraphics[width=1.1\textwidth]{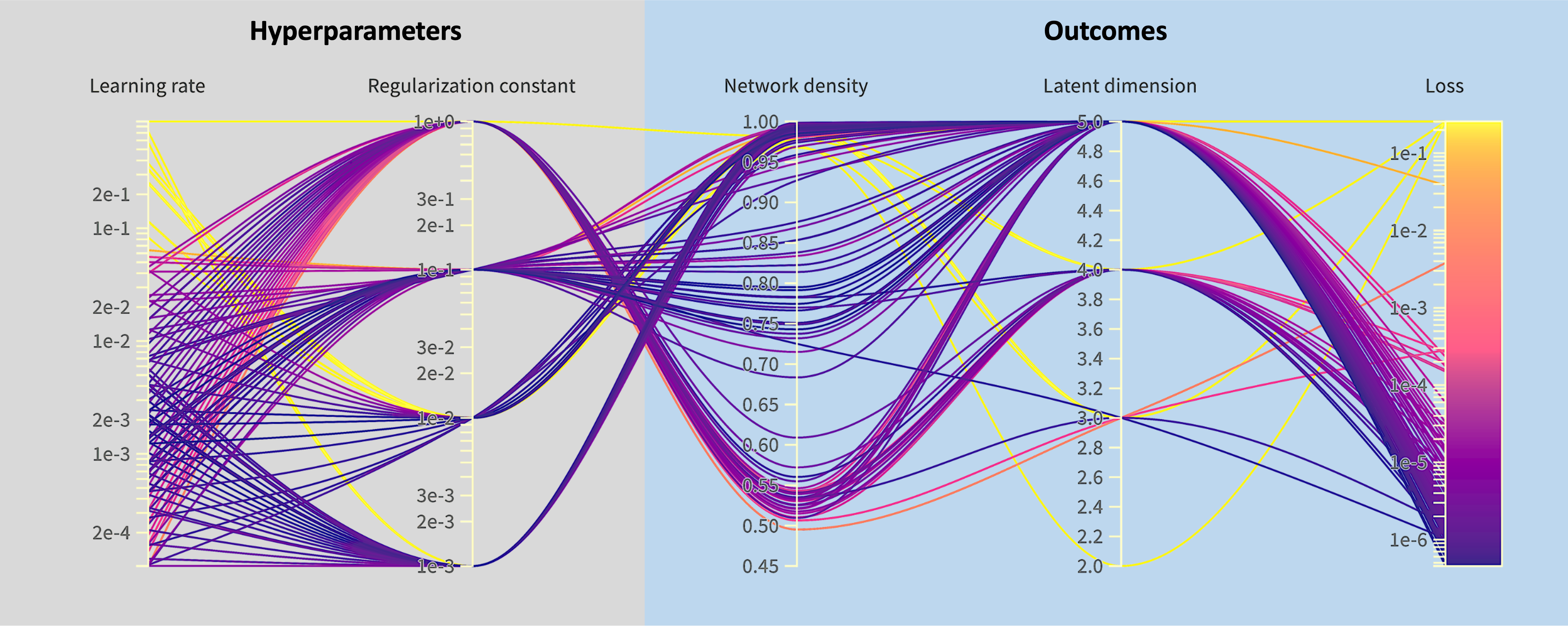}
    }
    \caption{Results for Wandb sweeps for the 1D diffusion equation and AdaBreg. Runs are started with the hyperparameters on the left, and the outcomes for the network with the lowest training loss achieved are shown on the right. One run with a large loss outlier for $\eta\approx 8.7$ and $\lambda=0.01$ has been omitted.}
    \label{fig:diffusion_adabreg_sweep}
\end{figure}

The last optimizer considered is AdaBreg. The four sweeps had $\lambda$ change one order of magnitude, similar to the LinBreg case. These four sweeps used the same lower bound of $\eta=\num{1e-4}$ for the learning rate. For the upper bound, the values $\eta=\num{10}$, $\eta=\num{0.1}$, $\eta=\num{10}$, and $\eta=\num{1}$ were used, respectively The results of these sweeps are visualized in \cref{fig:diffusion_adabreg_sweep}. The lowest loss is achieved for $\eta\approx\num{e-3}$ and $\lambda=0.01$. Many runs have a slightly worse but a similar loss. The run that stands out is with $\eta\approx\num{4e-3}$ and $\lambda=0.01$. It has a loss roughly a factor of $3$ from the best loss, but has a $40\%$ reduction in latent dimension and a $45\%$ reduction in the number of parameters. Hence, we take as learning rate $\eta=\num{4e-3}$ and as regularization constant $\lambda=1$.

\subsubsection*{Comparison}
To compare LinBreg and AdaBreg with SGD and Adam, we execute 10 runs with the chosen parameters and listed the performance of the run with the lowest testing loss in \cref{tab:1d-diffusion}. The adaptive methods have the lowest loss with factor of 2 between them. LinBreg and SGD are trailing behind also with a factor of 2 between them; this time with the factor favouring the Bregman-based method. Both Bregman methods use considerably fewer ($\approx78\%$ fewer) parameters than the classical methods. Only LinBreg uses fewer latent dimensions. AdaBreg outperforms Adam due to the fewer parameters used. Whether AdaBreg outperforms LinBreg as well depends on the preference for a 27-fold reduction in loss versus a 40\% reduction in latent dimensions.

\begin{table}[h]
    \centering
    \resizebox{\textwidth}{!}{%
    \Large
    \begin{tabular}{l|cccccccc}
         &
          \multicolumn{1}{c}{learning rate ($\eta$)} &
          \multicolumn{1}{c}{regularization ($\lambda$)} &
          \multicolumn{1}{c}{init density (\%)} &
          \multicolumn{1}{c}{\#params} &
          \multicolumn{1}{c}{latent dim} &
          \multicolumn{1}{c}{loss (training)} &
          \multicolumn{1}{c}{loss (testing)} \\ \hline\hline
        SGD         & $\num{5e-5}$ & - & 100 & $\num{12850}$ & 5 &  $\num{1.2e-4}$ &  $\num{1.2e-4}$ \\
        Adam        & $\num{1.5e-3}$ & - & 100  & $\num{12850}$ & 5 & $\num{1.1e-6}$ & $\num{1.1e-6}$ \\ \hdashline \noalign{\vskip 0.5ex}
        LinBreg     & $\num{e-3}$ & 1 & 20  & \textbf{$\num{2877}$}  & \textbf{3} & $\num{7.3e-5}$ & \textbf{$\num{6.0e-5}$} \\
        AdaBreg     & $\num{4e-3}$ & 1 &  20 &  \textbf{$\num{2849}$} &  \textbf{5} & $\num{2.1e-6}$ &   \textbf{$\num{2.2e-6}$}
    \end{tabular}%
    }
    \caption{Comparison in terms of losses, latent dimensions and number of parameters of the best models for the 1D diffusion equation using the four optimizers. POD achieves a testing loss of $\num{1.8e-8}$ with $5$ modes.}
    \label{tab:1d-diffusion}
\end{table}

\subsection{1D advection}\label{sec:advection}
The second example we consider is the 1D advection equation with periodic boundary conditions, which is given by 
\begin{subequations}\label{eq:advection_full}
    \begin{align}
         \partial_t u(x,t)+\mu_{\text{adv}}\partial_x u(x,t) &= 0 \\
         u(0,t) &= u(2,t) \\
         u(x,0) &= g(x)
    \end{align}
\end{subequations}
on the interval $I=[0,2)$ for $t\in [0,T]$ with $\mu_{\text{adv}}>0$, where 
\begin{equation}\label{eq:initial_condition_adv}
    g(x) = \frac{1}{\sqrt{0.002\pi}}e^{-\frac{(x-0.2)^2}{0.002}},
\end{equation}
i.e. a Gaussian with mean $0.2$ and variance $0.001$. The exact solution is given by
\begin{equation}\label{eq:advection_solution}
    u(x,t) = g((x-\mu_{\text{adv}}t) \mod 2),
\end{equation}
where the modulus is needed due to the periodic boundary conditions. Since this wave is travelling and preserves its shape, POD will require many modes to capture the solution accurately. When allowing for nonlinear maps, at most 2 modes are necessary. From a manifold perspective, an interval with periodic boundary conditions can be represented as a circle. Since the solution preserves its shape, we can map the solution $u(\cdot,t)$ to one point on this circle. A point on a circle requires at most two coordinates to describe. Hence, at most 2 coordinates are required to describe the function. When $T$ and $\mu_{\text{adv}}$ are chosen such that $u(2,t)=0$ for all $t\in [0,T]$, then we only need one coordinate to describe the solution. Moreover, the maps from the full-order solution to the one-dimensional latent space described by this one coordinate can be written down explicitly. Namely, the relations between FOM solutions $u(\cdot,t)$ and latent solution $z(t)$ are given by
\begin{subequations}\label{eq:advection_continuous_autoencoder}
\begin{align}
    \phi_{\text{enc}}&: u(\cdot,t) \mapsto z(t) = \int_\X x \,u(x,t)dx, \\
    \phi_{\text{dec}}&: z(t) \mapsto u(\cdot,t) = g(\cdot - z(t))
\end{align}
\end{subequations}
in which the univariate function $g$ is applied componentwise, since the expectation of a Gaussian is its mean. The corresponding latent dynamics satisfy
\begin{subequations}\label{eq:advection_latent}
\begin{align}
    \partial_t z(t) &= -\mu_{\text{adv}}, \\
    z(0) &= 0.2
\end{align}
\end{subequations}
with its exact solution given by
\begin{equation}
    z(t) = 0.2-\mu_{\text{adv}}t.
\end{equation}
The discrete versions of \eqref{eq:advection_continuous_autoencoder} are given by
\begin{subequations}\label{eq:constructed_autoencoder}
\begin{align}
    \bm{\phi}_{\text{enc}}&: \vb{u}_t \mapsto z_t=\frac{2}{N_x}\vb{x}^\intercal\vb{u}_{t},\\
    \bm{\phi}_{\text{dec}}&: z_t \mapsto g(\vb{x}-z_t).
\end{align}
\end{subequations}
The autoencoder $\bm{\phi}=\bm{\phi}_{\text{enc}}\circ \bm{\phi}_{\text{dec}}$ can be implemented using a no hidden layer neural network for $\bm{\phi}_{\text{enc}}$ and a three hidden layer neural network for $\bm{\phi}_{\text{dec}}$. Therefore, applying AdaBreg to this case will show how close the method can bring the encoder to the actual sparse latent space.
\begin{figure}
    \centering
    \includegraphics[width=0.48\textwidth]{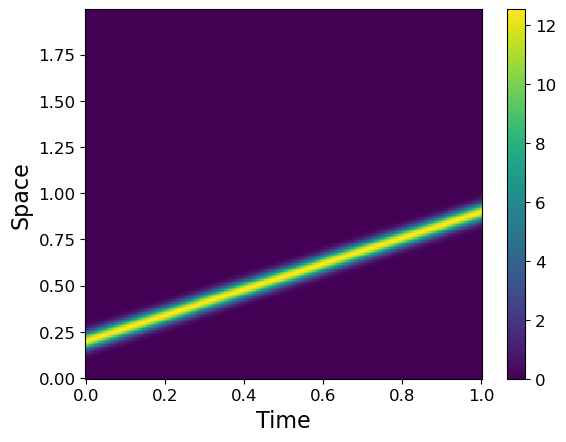}
    \includegraphics[width=0.48\textwidth]{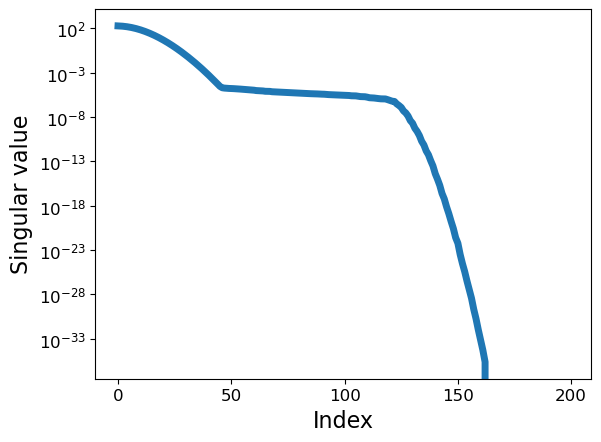}
    \caption{Numerical solution on the left for the diffusion equation with $\mu_{\text{adv}}=0.6$, and  on the right the singular value decay of the snapshot matrix corresponding to the shown numerical solution.}
    \label{fig:advection}
\end{figure}

For the FOM, the interval $I=[0,2)$ was divided into $N_x=256$ nodal points and $T=1$ taken with $N_t=200$ time steps. This implies the time step size $\Delta t$ and grid size $\Delta x$ are given by $\Delta t =1/N_t\approx\num{5.0e-3}$ and $\Delta x=1/N_x \approx\num{7.8e-3}$, respectively. For the advection parameter, we took $\mu_{\text{adv}}\in \Set{0.6,0.9,1.2}$ for the training set and $\mu_{\text{adv}}\in \Set{1.05}$ for the testing set. Combined with the chosen initial condition, this ensures that at no time $t\in [0,T]$ the wave touches the boundary. For the snapshots the exact solution \eqref{eq:advection_solution} was used, i.e.
\begin{equation}
    \vb{u}_n = g(\vb{x}+n\mu_{\text{adv}}\Delta t)
\end{equation}
where $g$ is as given in \eqref{eq:initial_condition_adv}. 

\subsubsection*{Architecture selection}
When POD is applied to the training set, $r=45$ singular values are required to get to a relative reconstruction error below $10^{-6}$. Since autoencoder allow for stronger compression, the bottleneck is set at $30$. Hence, autoencoders are set to have layers sizes $(256, 128, 50, 30, 50, 128, 256)$. 

For both training and testing, a batch size of 32 is used. Training is done for $1000$ epochs. 

\subsubsection*{Hyperparameter optimisation}
To find the best hyperparameters, four experiments have been conducted: one for each optimizer. These require one sweep for SGD as well as Adam and four sweeps for LinBreg as well as AdaBreg.

\begin{figure}[h!]
    \centering
    \begin{minipage}{0.48\textwidth}
        \includegraphics[width=\textwidth]{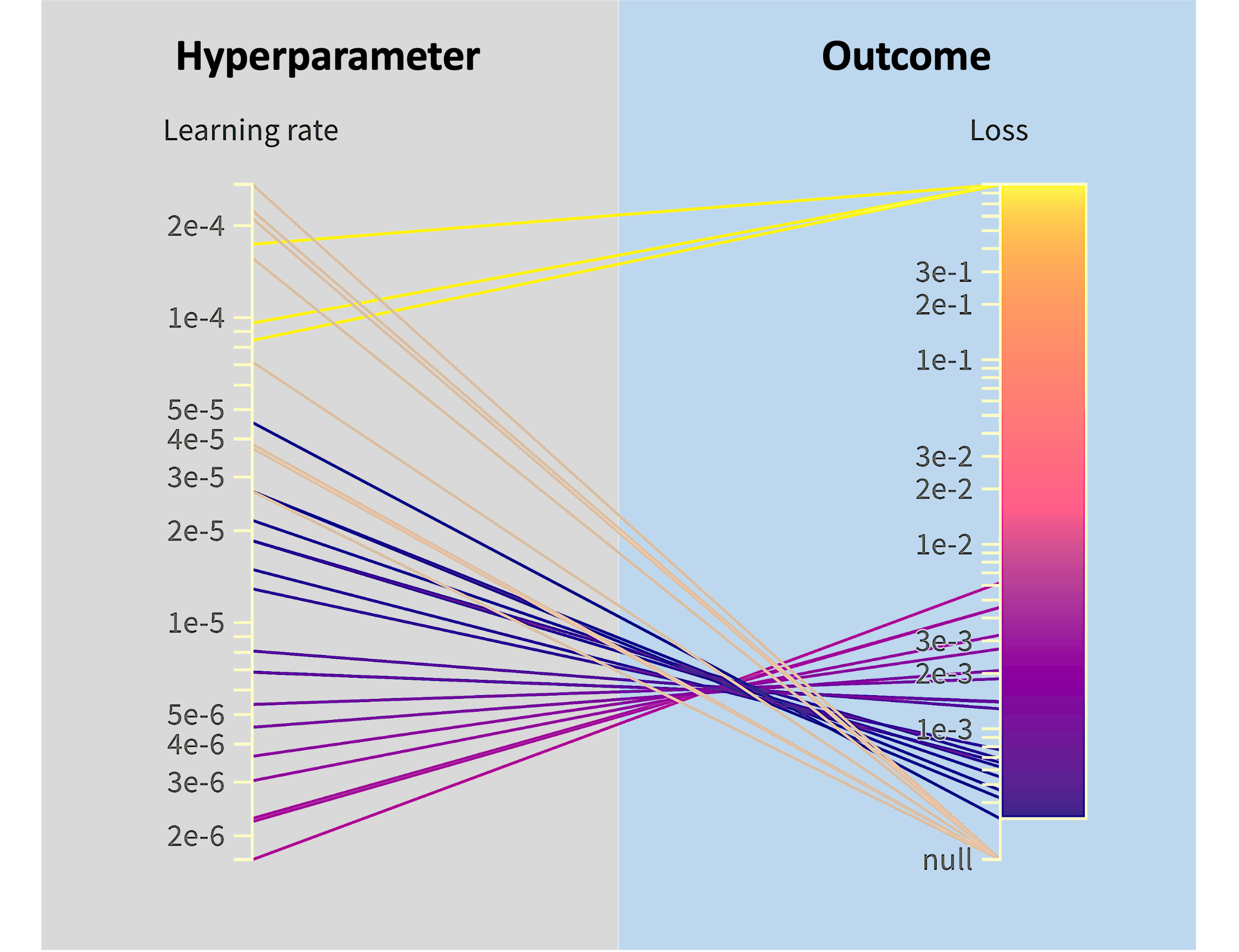}
        \caption{Wandb sweep for SGD. Runs are started with the learning rate on the left, and the lowest training loss achieved is shown on the right. The lowest loss is achieved for $\eta\approx\num{5e-5}$.}
        \label{fig:advection_sgd_sweep}
    \end{minipage}
    \begin{minipage}{0.48\textwidth}
        \includegraphics[width=\textwidth]{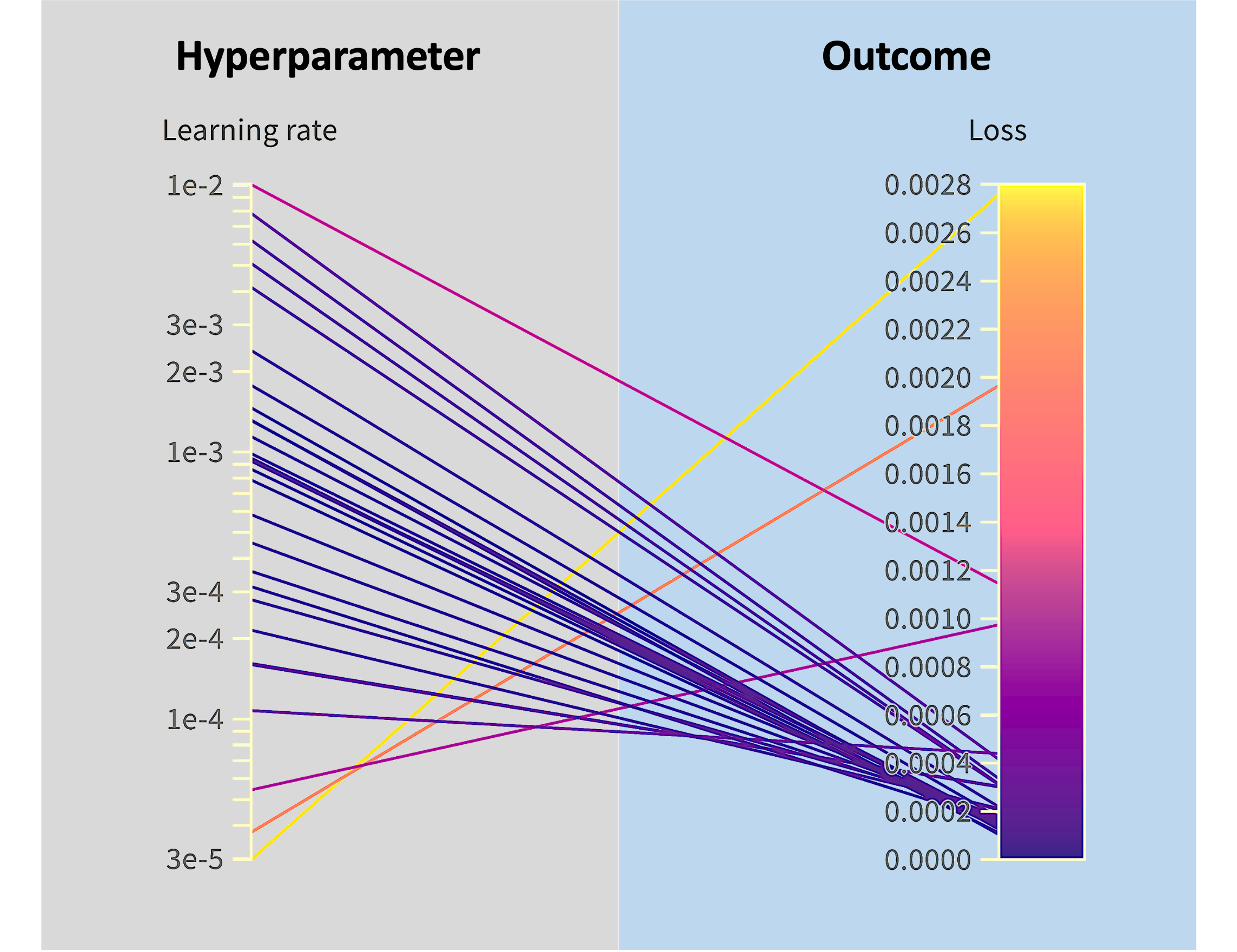}
        \caption{Wandb sweep for Adam. Runs are started with the learning rate on the left, and the lowest training loss achieved is shown on the right. The lowest loss is achieved for $\eta\approx\num{1.0e-3}$.}
        \label{fig:advection_adam_sweep}
    \end{minipage}
\end{figure}

For SGD, the prior runs from $\eta\approx \num{e-6}$ to $\num{e-4}$. The runs are shown in \cref{fig:advection_sgd_sweep}. From $\eta\approx\num{e-5}$ to $\eta\approx\num{5e-5}$, the runs all have similarly low losses around $\num{3e-4}$. Higher learning rates give unstable training and lower training rates give a higher loss. Hence, we choose $\eta=\num{4.5e-5}$.

The second optimizer considered is Adam. The bounds for the learning rate are from $\eta\approx \num{e-5}$ to $\num{e-2}$. From $\eta=\num{1e-5}$ to $\num{1e-3}$ the loss is decreasing, after which it starts increasing again. Hence, we take $\eta=\num{1e-3}$. 

\begin{figure}[h!]
    \centerline{
    \includegraphics[width=1.1\textwidth]{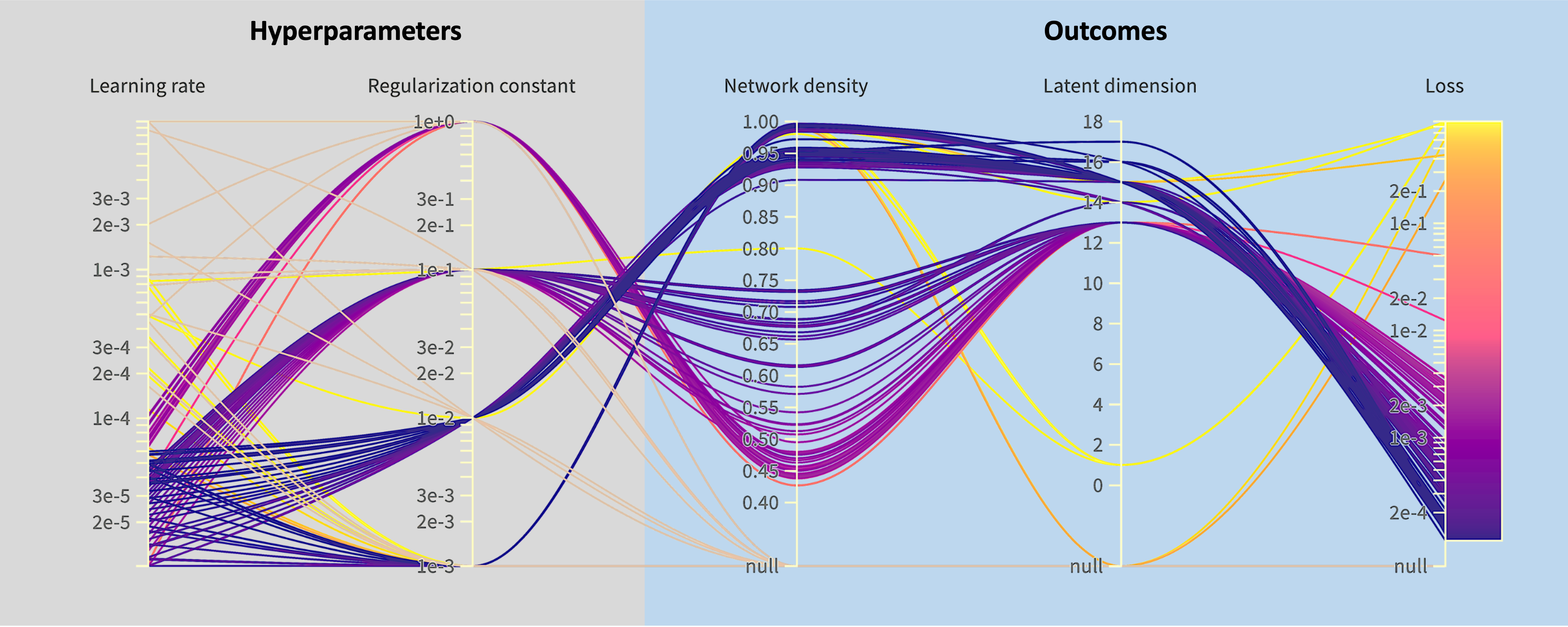}
}   
\caption{Results for Wandb sweep for the 1D advection equation and LinBreg. Runs are started with the hyperparameters on the left, and the outcomes for the network with the lowest training loss achieved are shown on the right.}
    \label{fig:advection_linbreg_sweep}
\end{figure}

Four different sweeps performed for LinBreg. For each of those sweeps, we changed the regularization constant $\lambda$. Just like with the 1D advection equation, $\lambda\in\Set{1, 0.1, 0.01, 0.001}$. The results of these sweeps are visualized in \cref{fig:advection_linbreg_sweep}. All runs use roughly half the allowed latent dimensions. The few runs that use significantly less do not perform well. The best performing network was for $\eta\approx\num{6e-5}$ and $\lambda=\num{e-2}$. This network is dense outside the latent dimension. The run with $\eta\approx\num{6e-5}$ and $\lambda=\num{e-2}$ has roughly $30\%$ fewer parameters. However, it has a loss a factor 6 higher. To stay competitive with the classical methods, we take the former hyperparameter set.

For the last optimizer AdaBreg, four sweeps were also performed. The same regularization constants were used as for LinBreg. In this case, there is a wide spread in the number of latent dimensions used. The run with $\lambda=0.01$ and $\eta\approx\num{1e-3}$ has the lowest loss of approximately $\num{9e-5}$. However, this run is eclipsed by the run with $\eta\approx\num{1.4e-3}$, $\lambda=0.1$ and a loss of approximately $\num{1.4e-4}$. This run has only 18 instead of 30 latent dimensions, and uses $30\%$ fewer parameters. Hence, we choose $\eta=\num{1.4e-3}$ and $\lambda=0.1$.

\begin{figure}[h!]
    \centerline{
    \includegraphics[width=1.1\textwidth]{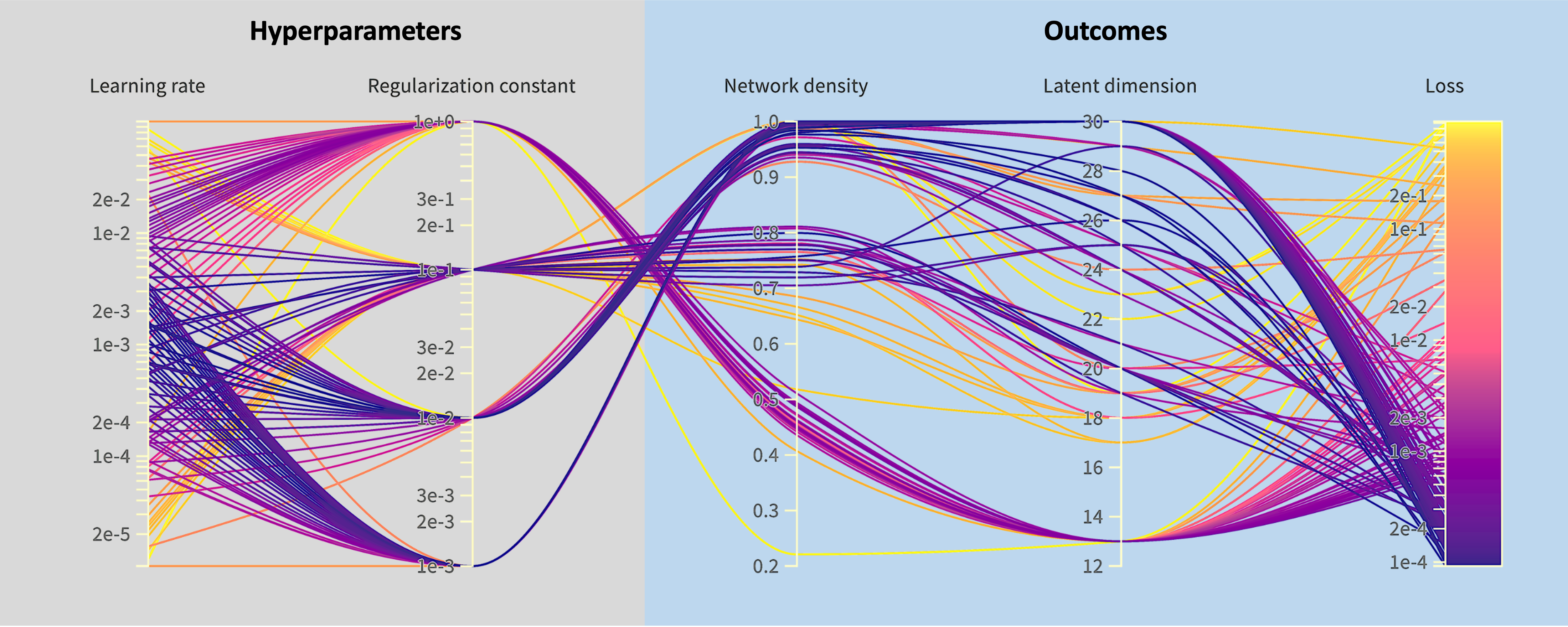}
}   
\caption{Results for Wandb sweep for the 1D advection equation and AdaBreg. Runs are started with the hyperparameters on the left, and the outcomes for the network with the lowest training loss achieved are shown on the right.}
    \label{fig:advection_adabreg_sweep}
\end{figure}

\begin{table}[h]
    \resizebox{\textwidth}{!}{%
    \begin{tabular}{l|cccccccc}
         &
          \multicolumn{1}{c}{learning rate ($\eta$)} &
          \multicolumn{1}{c}{regularization ($\lambda$)} &
          \multicolumn{1}{c}{initial density (\%)} &
          \multicolumn{1}{c}{\#parameters} &
          \multicolumn{1}{c}{latent dim} &
          \multicolumn{1}{c}{loss (training)} &
          \multicolumn{1}{c}{loss (testing)} \\ \hline\hline
        SGD         & $\num{4.5e-5}$ & - & 100 & $\num{85,760}$ & 30 & $\num{3.1e-4}$ & $\num{3.6e-4}$  \\
        Adam        & $\num{1.0e-3}$ & - & 100  & $\num{85,760}$ & 30 & $\num{1.6e-4}$ & $\num{1.7e-4}$  \\ \hdashline \noalign{\vskip 0.5ex}
        LinBreg     & $\num{6.0e-5}$ & $\num{0.01}$ & $20$ & $\num{75,008}$ & 17 & $\num{1.1e-4}$ & $\num{1.3e-4}$ \\
        AdaBreg     & $\num{1.4e-3}$ & $\num{0.1}$ &  $\num{20}$ & \textbf{$\num{52,398}$} & 13 & \textbf{$\num{1.2e-4}$} & \textbf{$\num{1.7e-4}$}
    \end{tabular}%
    }
    \caption{Comparison in terms of losses, latent dimensions and number of parameters of the best models for the 1D advection equation using the four optimizers. POD with 45 modes achieves a testing loss of \num{9.2e-7}.}
    \label{tab:1d-advection}
\end{table}

\subsubsection*{Comparison}
To compare LinBreg and AdaBreg with SGD and Adam, we execute 10 runs with the chosen parameters and listed the performance of the run with the lowest testing loss in \cref{tab:1d-advection}. All four optimizer yields networks with similar test loss. But, the network produced with AdaBreg needs less than half of the latent dimensions and approximately $40\%$ fewer parameters. Meaning that AdaBreg produced the best model.

\subsection{2D reaction-diffusion}
\begin{figure}[h!]
    \centering
    \begin{minipage}{0.48\textwidth}
        \includegraphics[width=\textwidth]{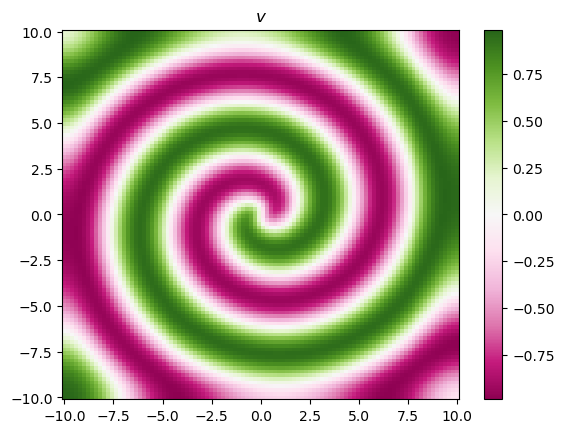}
        \caption{Solution $\vb{v}$ of the 2D reaction-diffusion equation at $t=5$.}
        \label{fig:reaction_diffusion_solution_v}
    \end{minipage}
    \begin{minipage}{0.48\textwidth}
        \includegraphics[width=\textwidth]{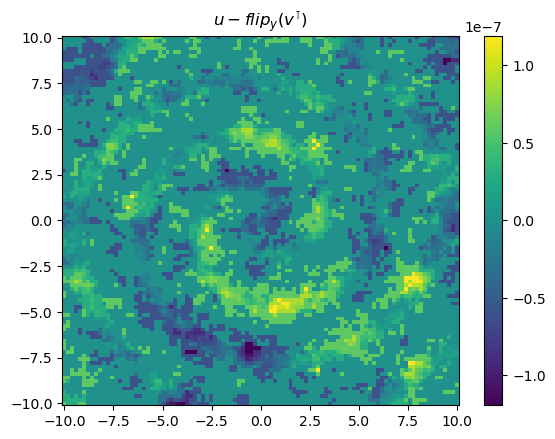}
        \caption{The solution $\vb{u}$ compared with $\text{flip}_y(\vb{v}^\intercal)$ at $t=5$, the latter being the transpose of $\vb{v}$ with its columns in reverse order.}
        \label{fig:reaction_diffusion_diff}
    \end{minipage}
\end{figure}

The third and last example we consider is a 2D reaction-diffusion problem whose governing equations are given by
\begin{subequations}
\begin{align}
    \partial_t u &= \mu_{diff,u}\Delta u + (1-u^2v^2) u + \beta u^2v^2 v \\
    \partial_t v &= \mu_{diff,v}\Delta v + (1-u^2v^2) v - \beta u^2v^2 u \\
    \grad v \cdot n &= \grad u \cdot n = 0 \label{eq:react_diffusion_bcs}\\
    u((x,y),0) &= \tanh{(\sqrt{x^2+y^2}\cos{(\angle(x+\i y)}-\sqrt{x^2+y^2}))} \\
    v((x,y),0) &= \tanh{(\sqrt{x^2+y^2}\sin{(\angle(x+\i y)}-\sqrt{x^2+y^2}))}
\end{align}    
\end{subequations}
on the domain $[-10,10]\cross[-10,10]$ for $t\in [T_1,T_2]$ with $\beta=1$. The full-order numerical solution is a spiral wave with the wave fronts getting smaller as time progresses. We have $T_1>0$, since the numerical solution takes some time to arrive at the spiral wave. If $T_2-T_1$ is small enough, then we can consider the wave fronts stable in size. In that case, we have a problem similar to the 1D advection in \cref{sec:advection}. The solution for $v$ matches the solution for $u$ after a 90-degree rotation and subsequent column reversal, and the solution for $u$ has a one dimensional latent manifold. The former is illustrated in \cref{fig:reaction_diffusion_diff}. There are no clear explicit maps for mapping the full-order solution to a latent solution and vice versa. Therefore, applying AdaBreg to this case will show how close the method can bring the encoder to the actual sparse latent space if the map is so complex.

For the full-order model, we took the diffusion parameters $\mu_{\text{diff},u}=\mu_{\text{diff},v}=1$, used a tensor grid over the domain $I=[-10,10]\cross[-10,10]$ with $N_x=100$ and $N_y=100$ nodal points respectively along the $x$- and $y$-directions, and set $T_2=1$ with $N_t=\num{50000}$ time steps. This implies the step size $\Delta t$ is set to $\Delta t\approx \num{2e-5}$. To evolve the wave in time from $t=0$ to $T_2$, we used operator splitting in a finite difference scheme. In particular,
\begin{subequations}
\begin{align}
    \vb{z}_{n} &= (\vb{u}_n\odot \vb{u}_{n})\odot(\vb{v}_{n}\odot \vb{v}_{n}) \\
    \vb{u}_{n+1} &= \vb{u}_{n} + \Delta t \bigg(\mu_{\text{diff},u}\bm{\Delta}\vb{u}_{n}+(1-\vb{z}_{n})\odot\vb{u}_{n} + \beta \vb{z}_{n}\odot\vb{v}_{n} \bigg) \\
    \vb{v}_{n+1} &= \vb{v}_{n}  + \Delta t \bigg(\mu_{\text{diff},v}\bm{\Delta}\vb{v}_{n}+\beta \vb{z}_{n}\odot\vb{u}_{n} + (1-\vb{z}_{n})\odot\vb{v}_{n} \bigg)
\end{align}
\end{subequations}
for all $n$, where $\bm{\Delta}$ is the second-order finite difference discretization of the 2D-Laplacian with ghost cells for the boundary conditions listed in \cref{eq:react_diffusion_bcs} and $\odot$ is the Hadamard product, i.e. the vectors to the left and to the right are pointwise multiplicated. 

The training and testing sets are subsamples of the numerical solution. The transient behaviour ends at $T_1=0.1$ so the first 5000 samples are discarded. The $v$ part of the solution is discarded of the remaining points. Subsequently, these points are subsampled with an interval of $36$ samples. Of these $\num{1250}$ subsamples, $\num{750}$ have been assigned to the training set, and $250$ to the testing set.

\begin{figure}[h!]
    \centering
    \begin{minipage}{0.48\textwidth}
        \includegraphics[width=\textwidth]{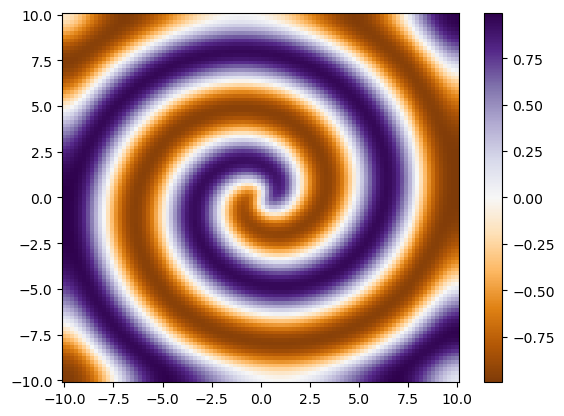}
        \caption{Solution $\vb{u}$ of the 2D reaction-diffusion equation at $t=5$.}
        \label{fig:reaction_diffusion_solution}
    \end{minipage}
    \begin{minipage}{0.48\textwidth}
        \includegraphics[width=\textwidth]{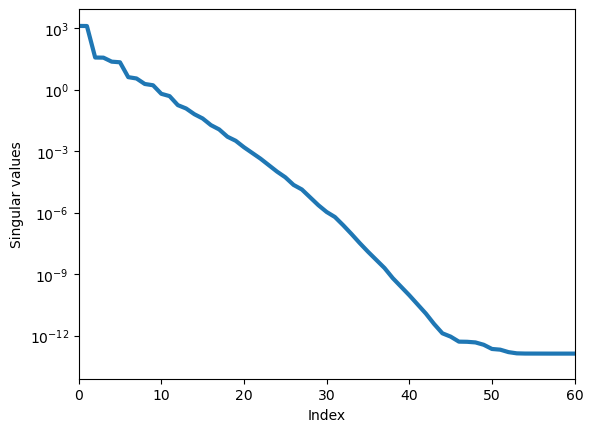}
        \caption{The first 60 singular values of the snapshot matrix with samples from the reaction-diffusion training set.}
        \label{fig:reaction_diffusion_train_set}
    \end{minipage}
\end{figure}

\subsubsection*{Architecture selection}
$r=10$ singular values are required to get a relative reconstruction error below $10^{-6}$ when using POD. The datasets only have the $u$ part, so the input dimensionality of the autoencoder should be \num{10000}. Hence, we use autoencoders with layers sizes $(\num{10000}, 200, 100, 10, 100, 200, \num{10000})$. 

\subsubsection*{Hyperparameter optimalisation}
To find the best hyperparameters, four experiments have been conducted: one for each optimizer. These require one sweep for SGD as well as Adam and four sweeps for LinBreg as well as AdaBreg.

\begin{figure}[h!]
    \centering
    \begin{minipage}{0.48\textwidth}
        \includegraphics[width=\textwidth]{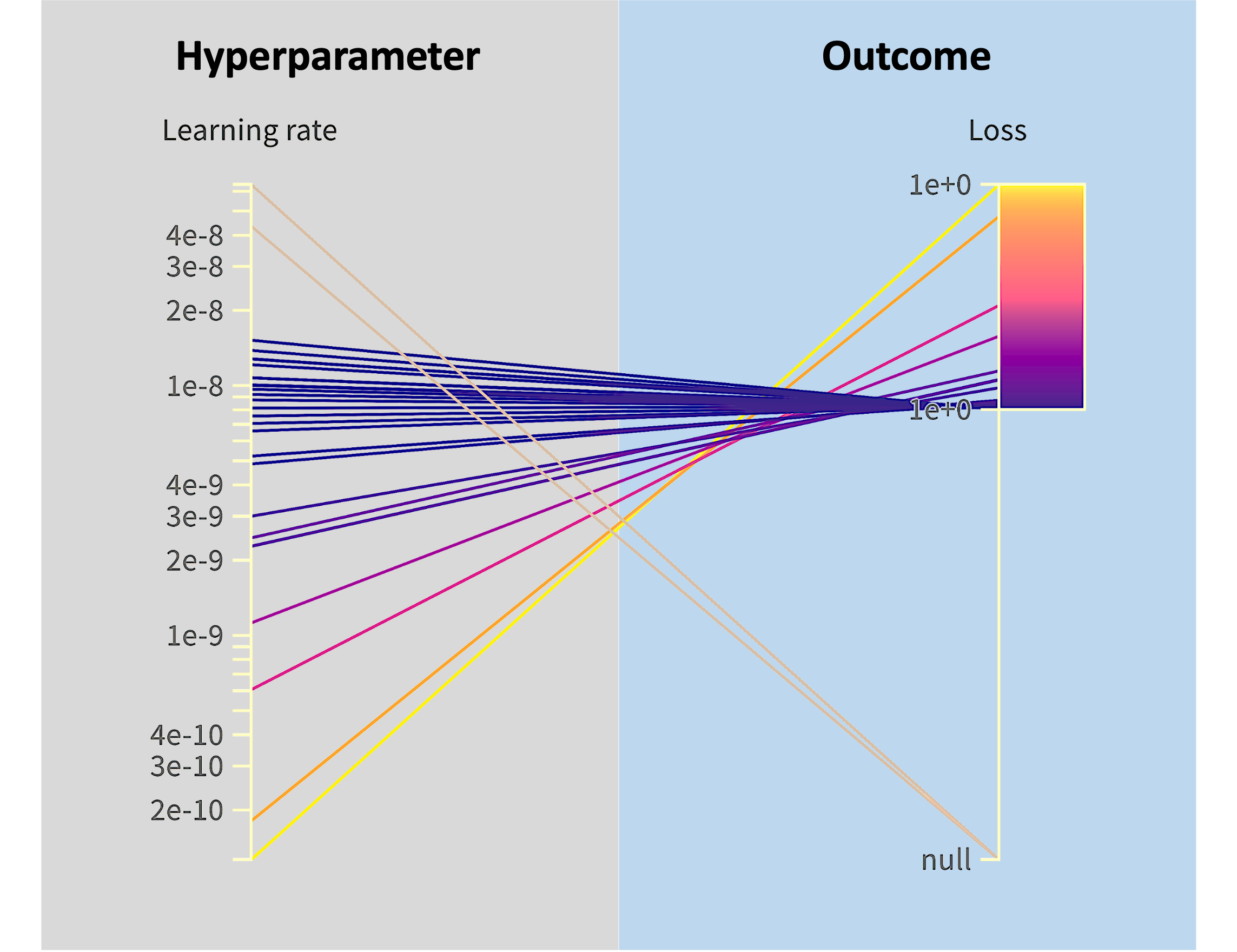}
        \caption{Wandb sweep for SGD. Runs are started with the learning rate on the left, and the lowest training loss achieved is shown on the right. The lowest loss is achieved for $\eta\approx\num{1e-8}$.}
        \label{fig:reaction_sgd_sweep}
    \end{minipage}
    \begin{minipage}{0.48\textwidth}
        \includegraphics[width=\textwidth]{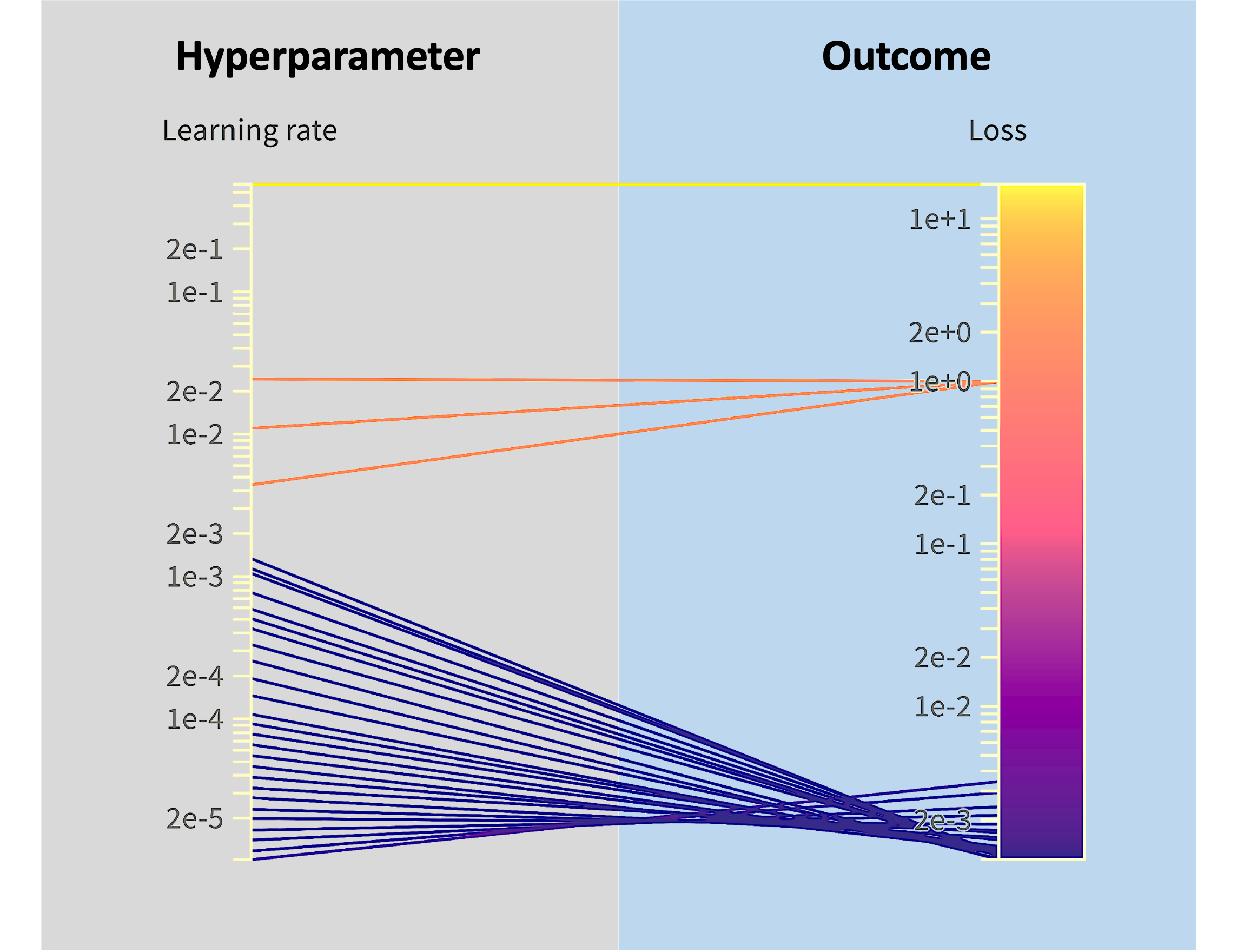}
        \caption{Wandb sweep for Adam. Runs are started with the learning rate on the left, and the lowest training loss achieved is shown on the right. The lowest loss is achieved for $\eta\approx\num{2e-3}$.}
        \label{fig:reaction_adam_sweep}
    \end{minipage}
\end{figure}

The first optimzer considered is SGD. This sweep has bound $\eta\approx \num{e-10}$ to $\num{5e-8}$. The runs are shown in \cref{fig:reaction_sgd_sweep}. From $\eta\approx\num{2e-10}$ the loss decreases, but all these runs have a loss of approximately $1$. When increasing the loss above $\num{4e-8}$, all runs have exploding gradients. This shows that with our current setup SGD fails to learn the data. Regardless, we choose for the final test $\eta=\num{4e-8}$.

The bounds for Adam are from $\eta\approx \num{e-5}$ to $\num{5e-1}$. The runs are shown in \cref{fig:reaction_adam_sweep}. From $\eta=\num{2e-5}$ to $\num{1.5e-3}$ the loss decreases. From $\eta\approx\num{4e-3}$ onward, the network fails to learn the data. Hence, we take $\eta=\num{1.5e-3}$. 

\begin{figure}[h!]
    \centerline{
    \includegraphics[width=1.1\textwidth]{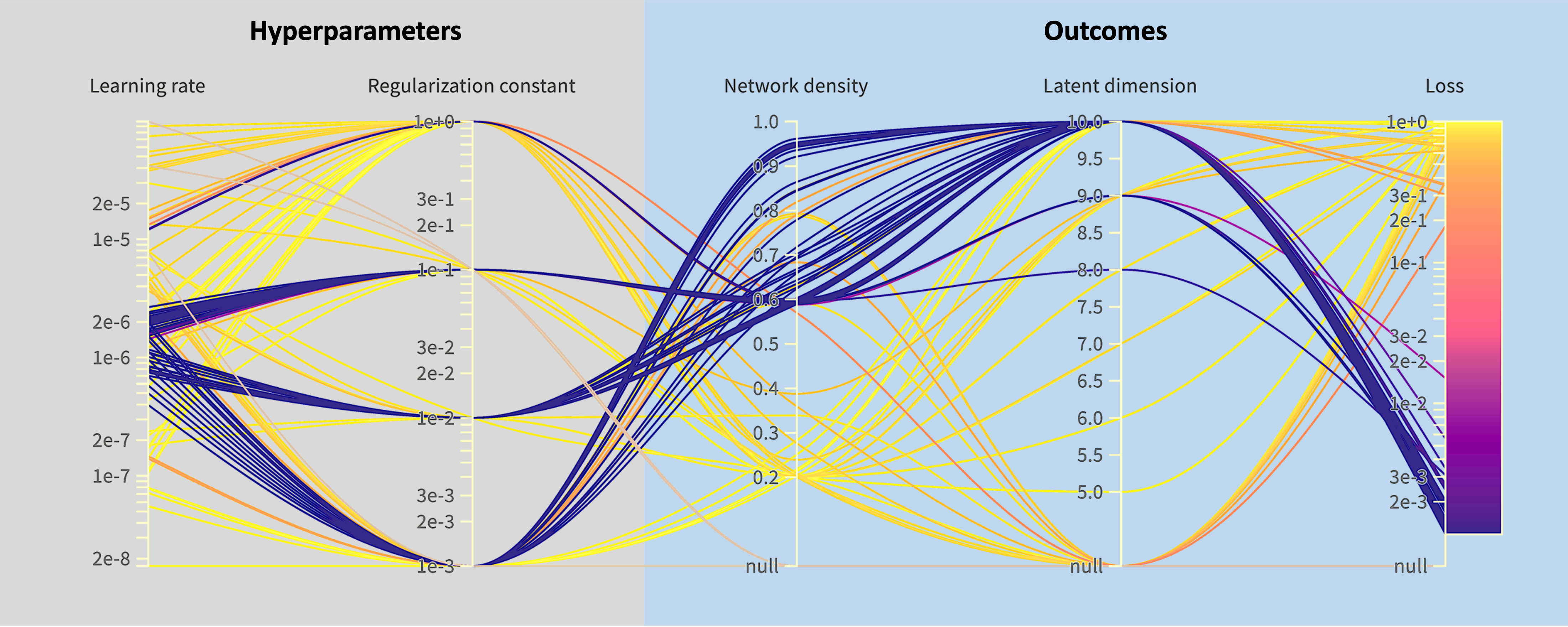}
}   
\caption{Results for Wandb sweep for the 2D reaction-diffusion equation and LinBreg. Runs are started with the hyperparameters on the left, and the outcomes for the network with the lowest training loss achieved are shown on the right.}
    \label{fig:reaction-diffusion_linbreg_sweep}
\end{figure}

For the third optimizer we consider, LinBreg, four different sweeps with the by now regular four values of $\lambda$ have been performed. For $\lambda=0.1$, $\lambda=0.01$ and $\lambda=0.001$ many runs have losses between $\num{e-3}$ and $\num{2e-3}$. Although they have similar losses, they are not similarly sparse. We choose $\eta=\num{2e-6}$ and $\lambda=0.1$ since these are the sparsest among them.

\begin{figure}[h!]
    \centerline{
    \includegraphics[width=1.1\textwidth]{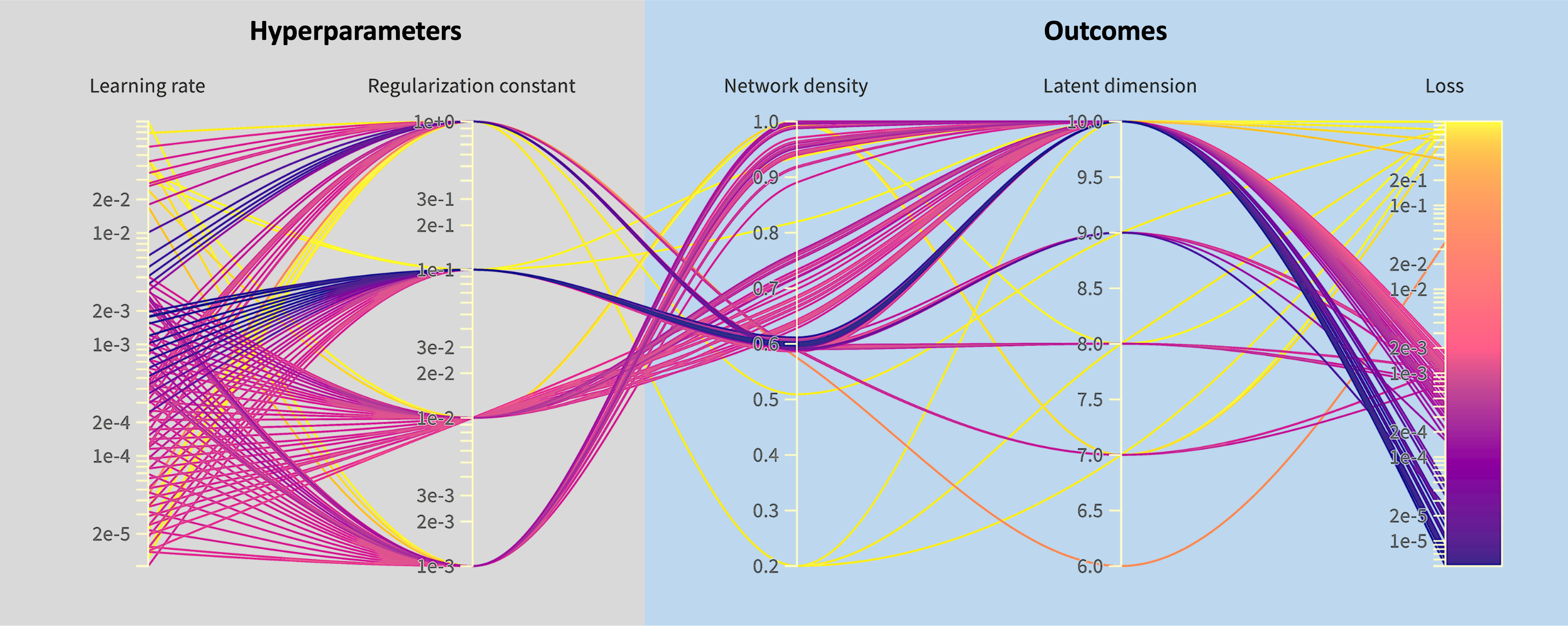}
}   
\caption{Results for Wandb sweep for the 2D reaction-diffusion equation and AdaBreg. Runs are started with the hyperparameters on the left, and the outcomes for the network with the lowest training loss achieved are shown on the right.}
    \label{fig:reaction-diffusion_adabreg_sweep}
\end{figure}

The results of our four sweeps for AdaBreg are shown in \cref{fig:reaction-diffusion_adabreg_sweep}. We see clearly that the best networks are those with $\lambda=0.1$ and $\eta\approx\num{2e-3}$. Hence, we take $\lambda=0.1$ and $\eta=\num{2e-3}$.

\subsubsection*{Comparison}
To compare LinBreg and AdaBreg with SGD and Adam, we execute 10 runs with the chosen parameters and listed the performance of the run with the lowest testing loss in \cref{tab:2d-reaction-diffusion}. SGD fails to properly with the chosen architecture. The remaining three optimizers yield networks with similar test losses. Both Bregman-based optimizer use significantly fewer parameters, with AdaBreg using half those of LinBreg. Meaning that is AdaBreg the best.

\begin{table}[h]
    \resizebox{\textwidth}{!}{%
    \begin{tabular}{l|cccccccc}
         &
          \multicolumn{1}{c}{learning rate ($\eta$)} &
          \multicolumn{1}{c}{regularization ($\lambda$)} &
          \multicolumn{1}{c}{init density (\%)} &
          \multicolumn{1}{c}{\#params} &
          \multicolumn{1}{c}{latent dim} &
          \multicolumn{1}{c}{loss (training)} &
          \multicolumn{1}{c}{loss (testing)} \\ \hline\hline
        SGD         & $\num{4e-8}$ & - & 100 & $\num{4042000}$ & 10 & $\num{1.0}$ & $\num{1.0}$ \\
        Adam        & $\num{1.5e-3}$ & - & 100  & $\num{4042000}$ & 10 & $\num{1.2e-3}$ & $\num{1.2e-3}$ \\ \hdashline \noalign{\vskip 0.5ex}
        LinBreg     & $\num{2e-6}$ & $\num{0.1}$ & 20  & $\num{1511652}$  & 9 & $\num{1.4e-3}$ & $\num{1.4e-3}$ \\
        AdaBreg     & $\num{2e-3}$ & $\num{0.1}$ &  20 & \textbf{$\num{769920}$} & \textbf{8} & $\num{1.1e-3}$ & \textbf{$\num{1.2e-3}$}
    \end{tabular}%
    }
    \caption{Comparison in terms of losses, latent dimensions and number of parameters of the best models for the 2D reaction-diffusion equation using the four optimizers. POD with 10 modes achieves a testing loss of \num{6.6e-8}.}
    \label{tab:2d-reaction-diffusion}
\end{table}

\section{Discussions}\label{sec:discussion}
In our presented algorithm, we use dense multilayer perceptrons instead of CNNs, choose a truncation level in the latent truncated-SVD, and chose a specific order of post-processing. In the following subsection, we will comment on those.

\subsection{Order of post-processing}
Our post-processing consists of two steps: latent truncated SVD and bias propagation. The order of these steps matters. With either order, all but the last layer of the encoder and first layer of the decoder will be the same. The difference is caused by zero rows in the matrix $\W^{L_{\text{enc}}-1}$. In general, $\W^{L_{\text{enc}}-1}$ will have zero rows. During the propagation, these rows and the corresponding columns from $\W^{L_{\text{enc}}}$ are removed during propagation. The spectral structure of $\W^{L_{\text{enc}}}$ is not preserved when columns are removed. This process is also called downdating [\cite{gu_downdating_1995}]. Empirically, it increases the latent space when applying POD on $\W^{L_{\text{enc}}}$. During the latent truncated SVD, the matrices $\W^{L_{\text{enc}}+1}$ and $\W^{L_{\text{enc}}+1}\vb{U}$ will have the same number of zero rows. Hence, doing the latent truncated SVD first preserves the structure for the propagation, whilst propagation destroys the structure for the latent truncated SVD.

\subsection{Computing latent truncated SVD error threshold}
When applying the truncated SVD in the latent space, we used a fixed error threshold during training. There was no guarantee that applying latent truncated SVD didn't make the network worse in reconstruction. We noted that the error is proportional to the Lipschitz constant of the decoder. Once we have a suitable estimate $\widehat{\text{Lip}}(\phiDec)$ of the Lipschitz constant of the decoder, we determine a suitable error criterion $\epsilon$ for the truncated-SVD from 
\begin{equation}
     \epsilon = \frac{c_{\text{tol}}\mathcal{L}(\theta^\dagger)}{\widehat{\text{Lip}}(\phiDec)}.
\end{equation}
In this, $\theta^\dagger$ collects the best parameters after training and $c_{\text{tol}}$ is a coefficient signifying how large the error due to the compression is allowed to be compared to the loss $\mathcal{L}(\theta^\dagger)$. 

Computing the Lipschitz constant of a neural network is NP-hard  [\cite{scaman_lipschitz_2019}]. There are several methods for estimating it. We will mention three. The first one is using the course upper bound
\begin{equation}\label{eq:latent_POD_upper_bound}
    \lip(\bm{\phi}_{\theta_{\text{post}},\text{dec}}) \leq \norm{V^\intercal}_2\lip(\sigma)^{L-L_{\text{enc}}}\prod_{\ell=L_{\text{enc}}}^L\norm{\W^\ell}_2 = \lip(\sigma)^{L-L_{\text{enc}}}\prod_{\ell=L_{\text{enc}}+1}^L s^\ell_0,
\end{equation}
where $\norm{\cdot}_2$ is the induced-vector norm and $s^\ell_0$ is the largest singular values of $\W^\ell$. The second one is using SeqLip from \cite{scaman_lipschitz_2019}. For the third one, we use that in many instances, the latent dimension is before the post-processing already relatively small. Hence, it is feasible to compute the Jacobians for (a subset of) the training data using forward-mode differentiation. The Lipschitz constant is bounded by the norm of the Jacobian. Hence, the constant can be bounded using
\begin{equation}
    \lip(\phiDec) \leq \max_{\bm{u}\in\X}\norm{\grad \phiDec(\bm{u})}_2.
\end{equation}

\subsection{Combining with CNN}
Our method enforces a particular sparsity pattern in feed-forward neural networks. Convolutional neural networks also have weight matrices with a particular pattern. However, these two patterns are incompatible with each other. Our method promotes zero rows in the matrices. Any zero row in a weight matrix of a convolutional neural network means the entire matrix must be zero due to their translation-equivariance. 

Although the sparsification techniques discussed in this paper cannot be applied directly, several options exist for applying Bregman iterations to CNNs. One way could be taking larger than usual kernels in the CNN. Then, an elementwise $L^0$ regularization can be used on the elements of the kernels to sparsify them. 

\section{Conclusions}
In this paper, we present a novel algorithm for finding sparse autoencoders. This algorithm is designed to find autoencoders that \emph{i}) are accurate in terms of reduction loss, \emph{ii}) have a small latent space, and \emph{iii}) are sparse. Specific loss functions and regularization were chosen such that the linearized Bregman iterations, upon which the algorithm is built, yield a network with these properties. To demonstrate the effectiveness of the method, we have this method on solution data generated from 3 different PDE problems: 1D-diffusion, 1D-advection, and 2D-reaction-diffusion. The best hyperparameters were found using Bayesian optimization and Wandb sweeps. The best resulting networks for each optimizer and PDE problems were compared. The summary of the comparison is shown in \cref{tab:summary}. The results show that our method can achieve a similar loss but with significantly fewer parameters.

\begin{table}[]
\hspace{-0.1\textwidth}\resizebox{1.2\textwidth}{!}{%
\begin{tabular}{l|lll|lll|lll}
        & \multicolumn{3}{c}{1D-diffusion} & \multicolumn{3}{c}{1D-advection} & \multicolumn{3}{c}{2D-reaction-diffusion} \\ \hline
 &
  \multicolumn{1}{c}{ldim} &
  \multicolumn{1}{c}{\#params} &
  \multicolumn{1}{c}{loss} &
  \multicolumn{1}{c}{ldim} &
  \multicolumn{1}{c}{\#params} &
  \multicolumn{1}{c}{loss} &
  \multicolumn{1}{c}{ldim} &
  \multicolumn{1}{c}{\#params} &
  \multicolumn{1}{c}{loss} \\ \hline\hline  \noalign{\vskip 0.5ex}
SGD     &    5 (100\%)   & $\num{12850}$ (100\%)& $\num{1.2e-4}$ & 30 (100\%) &     $\num{85.760}$ (100\%)& $\num{5.5e-4}$      &  10 (100\%)   &  $\num{4042000}$ (100\%)  &  $\num{1.0}$  \\
Adam    &    5 (100\%)   &$\num{12850}$ (100\%)& $\num{1.1e-6}$ &  30 (100\%)    &    $\num{85.760}$ (100\%)& $\num{1.7e-4}$      &   10 (100\%)  &  $\num{4042000}$ (100\%)  &  $\num{1.2e-3}$ \\ \hdashline \noalign{\vskip 0.5ex}
LinBreg &    \textbf{3 (60\%)}   & \textbf{$\num{2877}$ (22\%)}& \textbf{$\num{6.0e-5}$} &  17 (57\%)  &       $\num{75.008}$ (87\%)& $\num{1.3e-4}$      &  9 (90\%)  &  $\num{1511652}$ (37\%)  & $\num{1.4e-3}$  \\
AdaBreg &    \textbf{5 (100\%)}   & \textbf{$\num{2849}$ (22\%)}& \textbf{$\num{2.2e-6}$} & \textbf{13 (43\%)} &       \textbf{$\num{52.398}$ (61\%)} & \textbf{$\num{1.7e-4}$}     &   \textbf{8 (80\%)}  & \textbf{$\num{769920}$ (19\%)}  &   \textbf{$\num{1.2e-3}$}
\end{tabular}%
}
\caption{Summary of the numerical results from tables \ref{tab:1d-diffusion}, \ref{tab:1d-advection} and \ref{tab:2d-reaction-diffusion} for the 1D diffusion, the 1D advection and the 2D reaction-diffusion equation respectively. ldim refers to the latent dimension of the autoencoder, parameters refers to the number of non-zero weights and biases and loss refers to the testing loss. POD achieves a relative testing loss of \num{1.8e-6} with 5 modes, \num{9.2e-7} with 45 modes and \num{6.6e-8} with $10$ modes respectively. For all three equations, LinBreg and AdaBreg have a similar loss as SGD and Adam, but have considerably fewer parameters and smaller latent dimension.}
\label{tab:summary}
\end{table}

\section*{Acknowledgement}
This research was in part supported by Sectorplan Bèta (the Netherlands) under the focus area \emph{Mathematics of Computational Science}. C.~B. and M.~G. also acknowledge support from the 4TU Applied Mathematics Institute for the Strategic Research Initiative \emph{Bridging Numerical Analysis and Machine Learning}. 

\printbibliography

\end{document}